\documentclass{amsart}
\usepackage{amsfonts}
\usepackage{amssymb, amstext, amscd, amsmath}
\usepackage{amsthm}
\usepackage{times}
\usepackage{txfonts}
\usepackage{indentfirst}
\usepackage[all]{xy}
\usepackage{subfigure}
\usepackage{graphicx}
\usepackage{tikz}
\usepackage{pifont}
\usepackage{bbm}
\usepackage{amsfonts}
\usepackage{mathrsfs}
\usepackage{lscape}
\usepackage{setspace}
\usepackage{graphicx}
\usepackage{subfigure}
\usepackage{color}
\usepackage{indentfirst}
\usepackage{url}
\usepackage[utf8]{inputenc}
\newtheoremstyle{noparens}%
{}{}%
{\itshape}{}%
{\bfseries}{.}%
{ }%
{\thmname{#1}\thmnumber{ #2}\mdseries\thmnote{ #3}}
\theoremstyle{noparens}
\allowdisplaybreaks

\usepackage[colorlinks = true,
			linkcolor=blue,
			urlcolor=blue,
			citecolor=blue,
			anchorcolor=blue,
			backref=page]{hyperref}

\newtheorem{theorem}{Theorem}[section]
\newtheorem{corollary}{Corollary}[section]

\newtheorem{lemma}{Lemma}[section]
\newtheorem{remark}{Remark}[section]
\newtheorem{definition}{Definition}[section]
\newtheorem{example}{Example}[section]

\numberwithin{equation}{section}
\numberwithin{equation}{section}

\date{\today}

\begin{document}

\title[A three-variable transcendental invariant of planar knotoids via Gauss diagrams]{A three-variable transcendental invariant of planar knotoids via Gauss diagrams}
	
\author{Wandi Feng}
\address{School of Mathematical Sciences, Dalian University of Technology, Dalian 116024, P. R. China}
\email{fengwandi2000@163.com}
	
\author{Fengling Li$^*$}
\address{School of Mathematical Sciences, Dalian University of Technology, Dalian 116024, P. R. China}
\email{fenglingli@dlut.edu.cn}
	
\author{Andrei Vesnin}
\address{Sobolev Institute of Mathematics of the Siberian Branch of the Russian Academy of Sciences, Novosibirsk 630090, Russia}
\email{vesnin@math.nsc.ru}

\thanks{F.\,L. supported in part by grants (No. 12331003 and No. 12071051) of NSFC; A.\,V. supported by the State Task to the Sobolev Institute of Mathematics (No. FWNF-2022-0004)}
	
\begin{abstract}
As a generalization of the classical knots, knotoids are equivalence classes of immersions of the oriented unit interval  in a surface. In recent years, a variety of invariants of spherical and planar knotoids have been constructed as extensions of invariants of classical and virtual knots. In this paper we introduce  a three-variable transcendental invariant of planar knotoids which is defined over an index function of a Gauss diagram. We describe properties of this invariant and show that it is a Vassiliev invariant of order one. We also discuss the Gordian distance between planar knotoids and provide lower bounds on the Gordian distance of homotopic planar knotoids by using the transcendental invariant.
\end{abstract}
	
	\maketitle
	

\section{Introduction} \label{sec1}

The theory of knotoids was  introduced in 2012 by Turaev~\cite{Tu12}. Let $\Sigma$ be an oriented surface. Diagrams of knotoids are generic immersions of the oriented unit interval into $\Sigma$, together with the under/over crossing information at double points. Knotoids are defined as the equivalent class of knotoid diagrams under isotopies and Reidemeister moves. If  $\Sigma = S^2$ then the knotoid is said to be \emph{spherical}, and if $\Sigma = \mathbb{R}^2$ then the knotoid is said to be \emph{planar}. The theory of spherical knotoids extends the theory of classical knots and also provides a nice diagrammatic approach to classical knot theory. Basic notions of knotoids have been studied comprehensively  in~\cite{Tu12}, including the introduction of several invariants of knotoids and the monoid of knotoids with relations to $\theta$-graphs. A number of knot invariants were extended to knotoids.  In~\cite{GK17}  G\"{u}g\"{u}mc\"{u} and Kauffman introduced virtual knotoids, constructed some new invariants of classical and virtual knotoids, and provided applications of these polynomial invariants. A nice survey of knotoids, braidoids and their applications can be found in~\cite{GKL19}. Definitions of hyperbolicity for both spherical and planar knotoids have been  provided in~\cite{A24}.

In this paper we consider the set of planar knotoids. It was shown in~\cite{Tu12} that there is a surjective map from the set of planar knotoids to the set of spherical knotoids induced by the inclusion of $\mathbb R^2$  into $S^2$. However, the map is not injective since there are knotoid diagrams in $\mathbb R^2$ which represent a non-trivial knotoid while they represent the trivial knotoid in $S^2$. Thus, the theory of planar knotoids differs from the theory of spherical knotoids. In~\cite{KL19} Kodokostas and Lambropoulou studied a relation between planar knotoids and arcs in $\mathbb R^3$ considered up to a rail isotopy observed in~\cite{GK17}. It is shown  in~\cite{GGLDSK} that the technique of planar knotoids provides more refined information about the knotting of a protein than  other methods. A systematic classification of all planar knotoids with at most five crossings was given by Goundaroulis, Dorier and Stasiak in~\cite{GDS19}. In particular, it was shown that the number of prime planar knotoids with five crossings is between $944$ and $950$, where 6 pairs of knotoid diagrams were not distinguished. Then it was shown in~\cite{MV} that 2 of the 6 unresolved pairs of planar knotoids can be distinguished by quantum invariants.

Polynomial invariants are known to be useful for  distinguishing knotted objects. In recent years, many index-type polynomial invariants have been constructed for knots, virtual knots and knotoids. We will discuss some of them in Section~\ref{sec2}.
	
The aim of this paper is to introduce a new three-variable transcendental invariant $H_D(t,y,z)$ of a planar knotoid and to describe its properties and applications. 

Let $D$ be a planar knotoid diagram and $G(D)$ be the Gauss diagram of $D$. Denote the set of chords in $G(D)$ by $C(G(D))$. The function $H_D(t,y,z)$ is given below by the formula~(\ref{eqn:H}), that is
$$	
  H_D (t, y, z) = \sum_{\substack{c\in C(G(D))\\ n\in \mathbb{N}}} \operatorname{sgn} (c) \left ( t^{\operatorname{Ind}_{c}^{n}(z)} - 1  \right )y^n,
$$
where $\operatorname{Ind}_c^n(z)$ is the $n$-th index function given by formula~(\ref{eqn:ind}).
	
\begin{theorem}\label{thm-inv}
Let $D$ be a planar knotoid diagram, then the function $H_D(t,y,z)$ is an invariant of $D$.
\end{theorem}

In virtue of Theorem~\ref{thm-inv}, if $D$ is a diagram of a  knotoid $K$, we will use the notation $H_K(t,y,z)$ for $H_D(t,y,z)$.
	
\begin{remark} {\rm
The definition of $H_K(t,y,z)$ for planar knotoids can be extended to virtual knotoids, as the corresponding Gauss diagrams are not affected by the generalized Reidemeister moves and $\Omega_v$-moves, see~\cite{GK17}.
}
\end{remark}

The following theorem presents some properties of $H_D(t,y,z)$.

\begin{theorem} \label{thm-prop}
Let $D$ be a knotoid diagram.
\begin{itemize}
\item[(1)] If $-D$ is the inverse of $D$, then $H_{-D} (t, y, z) = H_{D} (t^{-1}, y, z)$.
\item[(2)] If $D^*$ is the mirror image of $D$, then $H_{D^*} (t, y, z) = -H_{D} (t^{-1}, y, z^{-1})$.
\item[(3)] If $D$ is zero height, then $H_D(t,y,z)=0$.
\end{itemize}
\end{theorem}
	
From point (2) of Theorem~\ref{thm-prop}, the following property holds immediately.

\begin{corollary}
If a planar knotoid diagram $D$ is invertible, then $H_D (t^{-1}, y, z) = H_D (t, y, z)$.
\end{corollary}

For point (3) of Theorem~\ref{thm-prop}, we recall that the height of a knotoid diagram is the least number of intersections between a diagram and an arc connecting its endpoints, where the minimum is taken over all representative diagrams and all such arcs are disjoint from crossings. A relation between the number of crossings  and the height of a knotoid is given by Korablev and Tarkaev in~\cite{KT21}. Note that the inverse of point (3) of Theorem~\ref{thm-prop} does not hold. That is, if $H_D(t,y,z)=0$, the planar knotoid diagram can have non-zero height, see Example~\ref{ex:1}.

The study of  Vassiliev invariants of knotoids was initiated in~\cite{MLK21}. It was shown in~\cite{MLK21} that for spherical knotoids there are non-trivial Vassiliev invariants of order one, in contrast to classical knot theory, where invariants of order one vanish. In the following theorem we discuss Vassiliev invariants for planar knotoids.

\begin{theorem} \label{thm-Vas}
Let $K$ be a planar knotoid, then $H_K(t,y,z)$ is a Vassiliev invariant of order one.
\end{theorem}
		
In~\cite{BG21} Barbensi and Goundaroulis considered $f$-distance, Gordian distance and band distance  of knotoids, corresponding to forbidden moves. We are interested in a Gordian distance related to crossing changes  for classical and virtual knots. Two planar knotoid diagrams are said to be homotopic, if they can be connected by a finite sequence of crossing changes and Reidemeister moves. For two homotopic planar knotoids $K$ and $K'$, we denote the Gordian distance between $K$ and $K'$ by $d_G (K, K')$. The following results allow to estimate $d_G (K, K')$ by comparing the invariants $H_K(t,y,z)$ and $H_{K'}(t,y,z)$.
	
\begin{theorem}\label{thm-dist}
If $K$ and $K'$ are two homotopic planar knotoids, then we have
$$
H_K(t,y,z)-H_{K'}(t,y,z) = \sum_{n\in\mathbb{N}}\left(\sum_{m\in\mathbb{N}}a_{n_m}\left( t^{ z_n^m} + t^{ -(z^{-1})_n^m}-2\right)\right)  y^n,
$$
where $a_{n_m} \in\mathbb{Z}$. Furthermore, $d_G(K,K')\geq\sum_{m\in\mathbb{N}}|a_{n_m}|$ for all $n\in\mathbb{N}$.
\end{theorem}
	
The paper is organized as follows. In Section~\ref{sec2}, we recall some basic concepts of knotoids such as Reidemeister moves, Gauss diagrams, Vassiliev invariants and Gordian distance. We then recall some index-type invariants of virtual knots and knotoids and define a three-variable transcendental function $H_D (t,y,z)$. In Section~\ref{sec3}, we give proofs of the Theorems~\ref{thm-inv}, \ref{thm-prop}, \ref{thm-Vas} and~\ref{thm-dist}. In Section~\ref{sec4}, we present several examples with calculations of $H_D(t,y,z)$ to illustrate properties of the invariant $H_D(t,y,z)$ of planar knotoids. 	
	
\section{Preliminaries} \label{sec2}

\subsection{Knotoids and Reidemeister moves} \label{subsec2.1}

We recall  the basic concepts of knotoids following~\cite{Tu12, GK17, GKL19}.
	
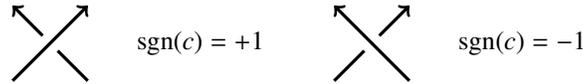
\begin{figure}[htbp]
	\begin{center}
		\begin{tikzpicture}
			\draw[black, very thick, ->] (0.4,0.6) -- (0,1);
			\node at (2.5,0.5) {$\operatorname{sgn} (c) = +1$};
			\draw[black, very thick] (0.6,0.4) -- (1,0);
			\draw[black, very thick, ->] (0,0) -- (1,1);
		\end{tikzpicture}
		\quad\quad
		\begin{tikzpicture}
			\draw[black, very thick, ->] (1,0) -- (0,1);
			\node at (2.5,0.5) {$\operatorname{sgn} (c) = -1$};
			\draw[black, very thick] (0,0) -- (0.4,0.4);
			\draw[black, very thick, ->] (0.6,0.6) -- (1,1);
		\end{tikzpicture}
	\end{center}
	\caption{The sign $\operatorname{sgn} (c)$ of a crossing $c$.} \label{fig1}
\end{figure}

Let $\Sigma$ be a surface. A \emph{knotoid diagram} $K$ in $\Sigma$ is a generic immersion of the interval $[0, 1]$ into the interior of $\Sigma$. Typically, knotoid diagrams are defined in $S^2$ or in $\mathbb R^2$, and the corresponding knotoid diagrams are called \emph{spherical} or \emph{planar}, respectively. In this paper we consider only the case of $\mathbb R^2$. Geometrically, a planar knotoid diagram is a curve in the plane that crosses itself transversely at each crossing point. The crossing points are finite and isolated, they are at most double points and have additional  information about over- or underpassing arcs. The sign of a crossing $c$ is defined as shown in Figure~\ref{fig1}.
We always assume that knotoid diagrams are \emph{oriented}, where the orientation is induced by the orientation of $[0,1]$ from $0$ to $1$. The images of $0$ and $1$ are called the \emph{tail} and the \emph{head} of the knotoid diagram, respectively, and both are different from the double points in the diagram.

\begin{figure}[htbp]
	\centering
	\vspace{8pt}
	\begin{tikzpicture}[scale=0.8]
		\draw[black, very thick] (0,1)--(0,3);
		\draw[black, very thick,<->] (0.2,2)--(0.8,2);
		\node[above] at (0.5,2) {$\Omega_1$};
		\draw[black, very thick] (1,3)--(1.5,2.1);
		\draw[black, very thick] (1.6,1.9)--(1.8,1.5);
		\draw[black, very thick] (1,1)--(1.8,2.5);
		\draw[black, very thick] (1.8,1.5) to [out=30,in=330] (1.8,2.5);	
	\end{tikzpicture}
	\qquad \qquad
	\begin{tikzpicture}[scale=0.8]
		\draw[black, very thick] (0,1)--(0,3);
		\draw[black, very thick] (0.5,1)--(0.5,3);
		\draw[black, very thick,<->] (0.7,2)--(1.3,2);
		\node[above] at (1,2) {$\Omega_2$};
		\draw[black, very thick] (1.5,1)--(2.5,2)--(1.5,3);
		\draw[black, very thick] (2.05,2.55)--(2.5,3);
		\draw[black, very thick] (1.95,1.55)--(1.5,2)--(1.95,2.45);
		\draw[black, very thick] (2.5,1)--(2.05,1.45);
	\end{tikzpicture}
	\qquad
	\qquad
\begin{tikzpicture}[scale=0.8]
	\draw[black, very thick] (0,0)--(0.9,0.9);
	\draw[black, very thick] (1.1,1.1)--(1.4,1.4);
	\draw[black, very thick] (1.6,1.6)--(2,2);
	\draw[black, very thick] (0,1.5)--(2,1.5);
	\draw[black, very thick] (0,2)--(0.4,1.6);
	\draw[black, very thick] (0.6,1.4)--(2,0);
	\draw[black, very thick,<->] (2.2,1)--(2.8,1);
	\node[above] at (2.5,1) {$\Omega_3$};
	\draw[black, very thick](3,0.5)--(5,0.5);
	\draw[black, very thick](3,0)--(3.4,0.4);
	\draw[black, very thick](3.6,0.6)--(3.9,0.9);
	\draw[black, very thick](4.1,1.1)--(5,2);
	\draw[black, very thick](3,2)--(4.4,0.6);
	\draw[black, very thick](4.6,0.4)--(5,0);
	\end{tikzpicture}
	\caption{Reidemeister moves $\Omega_1$, $\Omega_2$ and $\Omega_3$.}
	\label{fig2}
\end{figure}
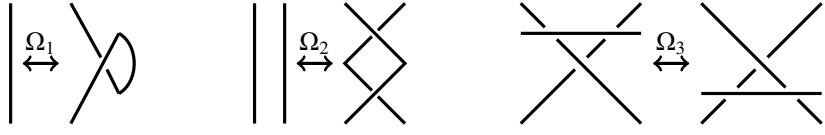

Consider the Reidemeister moves for knotoid diagrams as shown  in Figure~\ref{fig2}. When arcs in knotoid diagrams are considered to be oriented, the corresponding moves are called the \emph{oriented Reidemeister moves}.

\begin{definition}{\rm
Two oriented planar knotoid diagrams $D$ and $D'$ are said to be \emph{equivalent} if each can be obtained from the other after a finite sequence of performing the oriented Reidemeister moves away from the tail and  head  together with the local planar isotopy moves. The corresponding equivalence classes are called \emph{planar knotoids}.}
\end{definition}

\begin{figure}[htbp]
\begin{center}
\tikzset{every picture/.style={line width=1pt}}
\scalebox{1.0}{
\begin{tikzpicture}[x=0.56pt,y=0.56pt,yscale=-1,xscale=1]
	\draw[black, very thick]    (70,675) -- (150,675) ;
	\filldraw[black] (150,675) circle (1.5pt);
	\draw[black, very thick]    (110,630) -- (110,665) ;
	\draw[black, very thick]    (110,685) -- (110,720) ;
	\draw[black, very thick, <->]    (170,675) -- (210,675) ;
		\node[above] at (190,675) {$\Phi_+$};
    \draw[black, very thick]    (320,630) -- (320,720) ;
	\draw[black, very thick]    (240,675) -- (300,675) ;
	\filldraw[black] (300,675) circle (1.5pt);
	\draw[very thick, black, <->]    (360,675) -- (400,675) ;
		\node[above] at (380,675) {$\Phi_-$};
	\draw[black, very thick]    (440,675) -- (470,675) ;
	\draw[black, very thick]    (480,630) -- (480,720) ;
	\draw[black, very thick]    (490,675) -- (520,675) ;
	\filldraw[black] (520,675) circle (1.5pt);
\end{tikzpicture}
}
\caption{Forbidden knotoid moves $\Phi_+$ and $\Phi_-$.}  \label{fig3}
\end{center}
\end{figure}
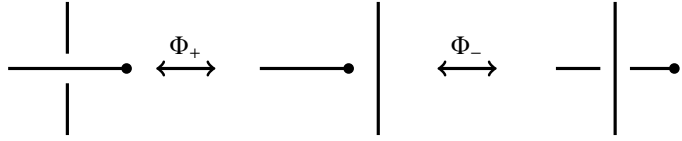

The moves consisting of pulling the arc adjacent to an endpoint over or under another arc, as presented in Figure~\ref{fig3} cannot be applied.

\begin{definition}{\rm
The moves $\Phi_+$ and $\Phi_-$ presented in Figure~\ref{fig3} are said to be \emph{forbidden}. }
\end{definition}
Notice that, if both $\Phi_+$ and $\Phi_-$ moves were allowed, any knotoid diagram in any surface could be clearly turned into the trivial knotoid diagram.
	
Planar knotoids can be split into two types: one is zero height planar knotoid, which has a diagram whose two endpoints are in the same region, and the other is nonzero height planar knotoid, which is a planar knotoid that does not have such diagram, see~\cite{Bat23}.

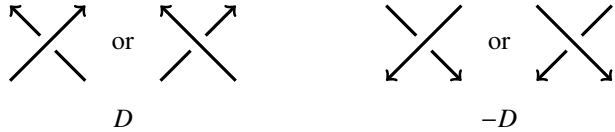
\begin{figure}[h]
\begin{center}
\begin{tikzpicture}
\draw[black, very thick, ->] (0.4,0.6) -- (0,1);
\node at (1.5,-0.5) {$D$};
\draw[black, very thick] (0.6,0.4) -- (1,0);
\draw[black, very thick, ->] (0,0) -- (1,1);
\node at (1.5,0.5) {or};
\draw[black, very thick] (2,0) -- (2.4,0.4);
\draw[black, very thick, ->] (3,0) -- (2,1);
\draw[black, very thick, ->] (2.6,0.6) -- (3,1);
\draw[black, very thick] (5,1) -- (5.4,0.6);
\draw[black, very thick, ->] (5.6,0.4) -- (6,0);
\draw[black, very thick, ->] (6,1) -- (5,0);
\node at (6.5,0.5) {or};
\node at (6.5,-0.5) {$-D$};
\draw[black, very thick] (7.6,0.6) -- (8,1);
\draw[black, very thick, ->] (7,1) -- (8,0);
\draw[black, very thick, ->] (7.4,0.4) -- (7,0);
\end{tikzpicture}
\end{center} \caption{Crossings of diagram $D$ and its inverse $-D$.} \label{fig4}
\end{figure}

For a knotoid diagram $D$, the diagram $-D$ obtained from $D$ by changing the orientation is called the \emph{inverse} of $D$. The corresponding changes of a crossing are shown in Figure \ref{fig4}. Observe that taking the inverse image does not change the sign of every crossing. It means that if $c$ is a crossing in $D$ and $c'$ is the image of $c$ in $-D$, then $\operatorname{sgn} (c') = \operatorname{sgn} (c)$.

\begin{figure}[h]
\begin{center}
\begin{tikzpicture}
\draw[black, very thick, ->] (0.4,0.6) -- (0,1);
\node at (1.5,-0.5) {$D$};
\draw[black, very thick] (0.6,0.4) -- (1,0);
\draw[black, very thick, ->] (0,0) -- (1,1);
\node at (1.5,0.5) {or};
\draw[black, very thick] (2,0) -- (2.4,0.4);
\draw[black, very thick, ->] (3,0) -- (2,1);
\draw[black, very thick, ->] (2.6,0.6) -- (3,1);
\draw[black, very thick] (5,0) -- (5.4,0.4);
\draw[black, very thick, ->] (5.6,0.6) -- (6,1);
\draw[black, very thick, ->] (6,0) -- (5,1);
\node at (6.5,0.5) {or};
\node at (6.5,-0.5) {$D^*$};
\draw[black, very thick] (7.6,0.4) -- (8,0);
\draw[black, very thick, ->] (7,0) -- (8,1);
\draw[black, very thick, ->] (7.4,0.6) -- (7,1);
\end{tikzpicture}
\end{center} \caption{Crossings of diagram $D$ and its mirror image $D^*$.} \label{fig5}
\end{figure}
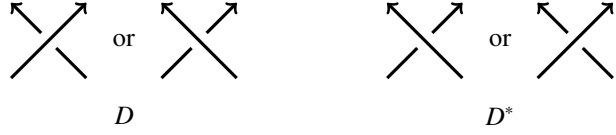

Let $D$ be a knotoid diagram, the diagram  $D^*$ obtained from $D$ by changing the over/under type of every crossing is called the \emph{mirror image} of $D$. The corresponding changes of a crossing are shown in Figure~\ref{fig5}. Observe that taking the mirror image does change the sign of every crossing to the opposite. It means that if $c$ is a crossing in $D$ and $c'$ is the image of $c$ in $D^*$, then $\operatorname{sgn} (c') = - \operatorname{sgn} (c)$.

\subsection{Index-type invariants of virtual knots and knotoids.} \label{subsec2.2}
In recent years, many index-type invariants have been constructed for virtual knots and knotoids. In a certain sense, these invariants are defined by taking the sign $\operatorname{sgn}(c)$ of a crossing $c$ and weight $W_K(c)$, associated with the crossings $c$ of a diagram of an oriented virtual knot or knotoid $K$, and the invariants take the form of a polynomial defined by the formula  
\begin{equation}
P_K (t) = \sum_{c} \operatorname{sgn} (c) \left ( t^{\displaystyle W_K (c) } - 1 \right ). \label{eqn:index}
\end{equation}
 Some authors refer to weights as indices. Among the index-type polynomial invariants of virtual knots we point out some invariants constructed by associating different weights.  In \cite{H10}, Henrich used  the intersection index $i(d)$ of an ordered virtual link associated with a smoothed crossing $d$ in a virtual knot diagram.  The odd writhe polynomial was constructed by Cheng~\cite{C14} by assigning some integers to arcs of Gauss diagram and taking the sum over odd crossings.  Cheng and Gao~\cite{CG13} defined the writhe polynomial using the Gauss diagram of virtual knots. Kauffman~\cite{K13} defined an affine index polynomial in terms of an integer labelling (known as Cheng coloring) of the virtual knot diagram. In~\cite{KPV} two-variable $n$-th $L$-polynomials and $F$-polynomials of oriented virtual knots were introduced by Kaur, Prabhakar and Vesnin. These polynomials arise from flat virtual knot invariants known as index value and $n$-th dwrithe. $F$-polynomials were tabulated by Ivanov and Vesnin in~\cite{IV} for oriented virtual knots with at most four classical crossings in a diagram. In \cite{HL}, Hao and Li constructed a family of polynomials $H_D^n(t,h,l)$ of virtual knots by some smoothing rules. In \cite{GIPV}, Gill, Ivanov, Prabhakar and Vesnin introduced weight functions for ordered orientable virtual and flat virtual links.  By considering three types of smoothing in classical crossings of a virtual link diagram and using appropriate weight functions, they provided a recurrent construction for new invariants.

Recently, based of the writhe polynomial and the affine index polynomial for virtual knots, Jeong~\cite{Je23} defined a virtual knot polynomial invariant $H_D(x,y)$ with two variables as follows
$$
H_D (x, y) = \sum_{\substack{c \in C (D) \\ n\in \mathbb{N}}} \operatorname{sgn} (c) \left ( x^{d_n (c)} -1 \right) y^n,
$$	
where $d(c)$ is a degree of crossing  $c$ of a virtual knot diagram $D$, and $d_n(c)$ is obtained by counting the number of different chords $c'$ transversing $c$ with $\gcd(d(c), d(c')) = n$. 

In~\cite{C17} Cheng defined a virtual knot invariant using transcendental functions, which extends all polynomial invariants mentioned above and found its relation to  other index-type polynomial invariants. Let $K$ be a virtual knot diagram and $G(K)$ the corresponding Gauss diagram. Choose a chord $c$ in $G(K)$ and assign it an index Ind($c$) based on the direction and sign of the other chords passing through it. According to the direction of these chords, the chords are divided into two classes. The index function $g_c(s)$ is defined as follows 
$$
g_c (s) = \sum_{i=1}^{n} w (r_i ) s^{ \phi ( \operatorname{Ind} (r_i) )} - \sum_{i=1}^{m} w(l_i) s^{ \phi ( - \operatorname{Ind} (\ell_i) )},
$$
and a virtual knot invariant $F_K(t,s)$ is given by the formula
$$
F_K (t, s) = \sum_{c_i} w (c_i) t^{\displaystyle g_{c_i} (s)} - w ( K).
$$
	
Correspondingly, based on virtual knot index type invariants, similar index type invariants are also defined on knotoids. In~\cite{GK17}  G\"{u}g\"{u}mc\"{u} and Kauffman defined the odd writhe for both classical and virtual knotoid diagrams as the sum of the signs of the odd crossings $D$. They also gave the weight $W_K(c)$ at the crossing $c$ and obtained the affine index polynomial which looks similar to~(\ref{eqn:index}) of a virtual or classical knotoid diagram~$K$.	
	
In~\cite{KIL18} Kim, Im and Lee defined the index polynomial for knotoids and the $n$-th polynomial invariants for virtual knotoids using Gauss diagram and provided some properties. The index polynomial $F_{G(D)}(t)$ is defined by
$$
F_{G(D)}(t)=\sum_{c\in C\left ( G\left ( D \right )  \right ) } \operatorname{sgn} (c) \left ( t^{\displaystyle i (c)} - 1 \right ),
$$
and the $n$-th polynomial invariants $Z_{G(D)}^{n} (t)$ are defined by
$$
Z_{G(D)}^{n} (t ) = \sum_{c \in C_n ( G (D) ) }\operatorname{sgn} (c) \left ( t^{\displaystyle d_n (c)} - 1 \right ).
$$
Feng and Li~\cite{FL23} used Gauss diagrams to define new two-variable invariant of knotoids, which is given  by
$$
F_D(u,v)=\sum_{c\in C(G(D))}\mathrm{sign}(c)(u^{g_c(v)}-1).
$$ 

In addition to index-type polynomial invariants, knotoids also have polynomial invariants in other forms as well. In~\cite{Bat22}, Bataineh introduced invariants for both spherical knotoids and planar knotoids through the one-to-one correspondence between knotoids and 2-polar knots. In~\cite{Bat23}, Bataineh, Batayneh and  Alkasasbeh introduced the winding signed sum polynomial of planar knotoids and gave an application to geometric invariants. 

\subsection{Gauss diagram and $H_D (t, y, z)$ function.}	
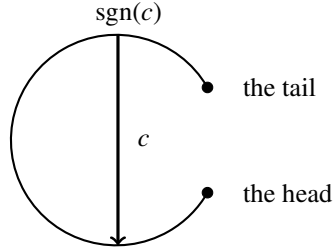
\begin{figure}[htbp] 
	\centering
	\begin{tikzpicture}[scale=0.7]
		\draw[black, thick] (0,0) arc (30:330:2);
		\filldraw[black] (0,0) circle (.1);
		\node at (0.5,0) [right]{the tail};
		\filldraw[black] (0,-2) circle (.1);
		\node at (0.5,-2) [right]{the head};
		\draw [black, very thick,->](-1.7,1) -- (-1.7,-3);
		\node [above] at(-1.5,1){sgn($c$)};
		\node[right] at (-1.5,-1){$c$};
	\end{tikzpicture}
	\caption{Elements of a Gauss diagram $G(D)$.}
	\label{fig6}
\end{figure}

The Gauss diagram $G(D)$ of a knotoid diagram $D$ is a counterclockwise  oriented arc with chords connecting preimages of crossings of $D$. The starting point of  $G(D)$ is called the \emph{tail} of $G(D)$, and the end point of $G(D)$ is called the \emph{head} of $G(D)$.  The chord of  $G(D)$ corresponding to a crossing $c \in D$ is also denoted by $c$ as well. A chord gets a direction from the overcrossing to the undercrossing. We also assign  a \emph{sign} $\operatorname{sgn} (c)$ to the starting point of a chord $c$ according to the sign of the corresponding crossing of $D$, see Figure~\ref{fig6}.

Note that we are considering chords whose endpoints belong to the arc between the tail and the head, but not to the whole circle. This is the difference between the Gauss diagram of a knotoid and the Gauss diagram of a virtual knot.	
	
\begin{figure}[htbp] 
	\centering
	\begin{tikzpicture}[scale=0.7]
		\draw[black, thick](0,0) arc (30:330:2);
		\filldraw[black] (0,0) circle (.1);
		\filldraw[black] (0,-2) circle (.1);
		\draw [black, very thick,<-](-1.7,1) -- (-1.7,-3);
		\draw [black, very thick,<-](-0.2,0.3) -- (-3.22,0.3);
		   \node at (-3.21,0.3) [left,scale=0.8]{$-$};
		   \node at (-3,0.6) [left,scale=0.8]{$+$};
		    \node at (1.2,0.3) [left,scale=0.8]{$r^-(c)$};
		    \node at (1,0.6) [left,scale=0.8]{$r^+(c)$};
		\draw [black, very thick,<-](-0.6,0.6) -- (-2.9,0.6);
		\node [below]at(-1.5,-3){sgn($c$)};
		\node[right] at (-1.5,-1){$c$};
		\draw [black, very thick,->](-0.2,-2.3) -- (-3.21,-2.3);
		\draw [black, very thick,->](-0.5,-2.6) -- (-2.9,-2.6);
		  \node at (.4,-2.3) [left,scale=0.8]{$+$};
		   \node at (.2,-2.6) [left,scale=0.8]{$-$};
		        \node at (-3.5,-2.2) [left,scale=0.8]{$\ell^+(c)$};
		         \node at (-3.3,-2.6) [left,scale=0.8]{$\ell^-(c)$};
	\end{tikzpicture}
	\caption{Positive and negative chords of $G(D)$ passing through a chord $c$.}
	\label{fig7}
\end{figure}
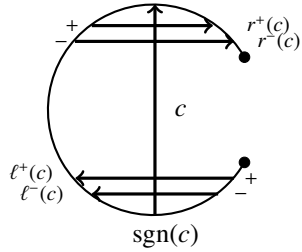

Let $C(G(D))$ be the set of chords in the Gauss diagram $G(D)$ of a knotoid diagram $D$. We define the \emph{degree} $d(c)$ of the chord $c \in C(G(D))$ by
\begin{equation}
d(c) = |r^+ (c)| -  |r^- (c)| - |\ell^+(c)| + |\ell^- (c)|, \label{eqn2.1}
\end{equation}
where $r^+(c)$ (or $r^-(c)$) is the set of chords with positive (or  negative) signs passing through $c$ from left to right along the direction of $c$, and $\ell^+(c)$ (or $\ell^-(c)$) is the set of chords with positive (or negative) sign passing through $c$ from right to left along the direction of $c$,  see Figure~\ref{fig7}. If chord $c$ is isolated, i.e. no other chord has non-empty intersection with $c$, then $d(c) = 0$. Note that the above defined degree $d(c)$ of a chord was used in~\cite{Je23}, also it was used in~\cite{C16} called the index of a chord. With $|s|$ we denote the cardinality of a set $s$.
	
For a fixed chord $c \in C(G(D))$ we define a function $\phi_c : \mathbb Z \to \mathbb Z_{|d(c)|}$ which corresponds to an integer $k$ its value by mod $|d(c)|$, namely,
\begin{equation}
\phi_c (k) = k \operatorname{mod} | d(c) |. \label{eqn2.100}
\end{equation}
with $\phi_c = \mathrm{id}$ if $d(c) = 0$. 
	
For a fixed chord $c \in C(G(D))$ we divide the chords that pass through the chord $c$ into two classes only according to their directions. Let $r(c)$ be the set of chords that pass through $c$ from left to right in respect to the direction of $c$, then  $r(c) = r^{+} (c) \cup r^{-} (c)$, and $\ell(c)$ be the set of chords that pass through $c$ from right to left with respect to the direction of $c$, then $\ell(c) = \ell^{+} (c) \cup \ell^{-} (c)$, see Figure~\ref{fig7}.

 For $c \in C(G(D))$ and $n \in \mathbb N$ we consider the set 
 $$
C_n(G(D), c) = \{ e \in C  ( G ( D )  ) \, \mid \, e \cap c \neq \emptyset, \,  \gcd ( d ( c ), d ( e) ) =n \},
$$
where $d(c)$ and $d(e)$ are degrees given by (\ref{eqn2.1}). Recall that the greatest common deviser $\gcd(a,b)$ of integers $a$ and $b$, at least one of which is nonzero, is the greatest positive integer which divides both; moreover, $\gcd(a,0) = \gcd(0,a) = |a|$ and $\gcd(0,0)=0$.

Let us split $C_n(G(D), c)$ into two subsets: one contains chords passing through $c$ from right to left,
$$
\ell^n(c) = \ell(c) \cap C_n(G(D), c),
$$
and another contains chords passing through $c$ from left to right,
$$
r^n(c) = r(c) \cap C_n(G(D), c).
$$

For  $c \in C(G(D))$ and $n \in \mathbb N$  we define the \emph{$n$-th index function} $\operatorname{Ind}_c^n(z)$ in variable $z$ by the formula 
\begin{equation}
\operatorname{Ind}_{c}^{n} (z) =  \sum_{e \in r^n(c)} \operatorname{sgn} (e) \, z^{\,  \phi_c (d(e)) } - \sum_{e \in \ell^n(c)} \operatorname{sgn} (e) \, z^{\,  \phi_c (-d(e))},   \label{eqn:ind}
\end{equation}
where $\phi_c (k)$ is given by (\ref{eqn2.100}). More precisely, if $d(c)=0$, then $\phi_c = \mathrm{id}$, hence
$$
\operatorname{Ind}_{c}^{n} (z) =  \sum_{e \in r^n(c)} \operatorname{sgn} (e) \, z^{d(e)} - \sum_{e \in \ell^n(c)} \operatorname{sgn} (e) \, z^{-d(e)}.
$$
If $d(c) = \pm 1$, then $\phi_c(k) = 0$ for any $k$, then 
$$
\operatorname{Ind}_{c}^{n} (z) =  \sum_{e \in r^n(c)} \operatorname{sgn} (e)  - \sum_{e \in \ell^n(c)} \operatorname{sgn} (e).
$$
 In general, $\operatorname{Ind}_c^n(z)$ takes values in  $\mathbb{Z} [ z^{\pm 1} ]/ ( z^{\mid d(c)\mid}-1 )$.

\begin{definition}
{\em
Let $D$ be a planar knotoid diagram and $G(D)$ its Gauss diagram. Denote by $C(G(D))$ the set of all chords in $G(D)$. Define a three-variable function $H_D(t, y, z)$ as follows}
\begin{equation}
H_D(t, y, z) = \sum_{\substack{c\in C(G(D))\\ n\in \mathbb{N}}} \operatorname{sgn} (c) \left( t^{  \operatorname{Ind}_{c}^{n} (z) } - 1 \right ) y^n. \label{eqn:H}
\end{equation}
\end{definition}

Note that if a knotoid diagram has a finite number of crossings, then the sum in the definition is finite since the number of pairs $(c,n)$ such that $\operatorname{Ind}_c^n (z) \neq 0$ is finite. 

\subsection{Singular planar knotoid diagrams.}

\begin{figure}[h]
\begin{center}
\begin{tikzpicture}
\draw[black, very thick, ->] (1,0) -- (0,1);
\draw[black, very thick, ->] (0,0) -- (1,1);
\filldraw [black] (0.5, 0.5) circle[radius=0.08];
\node at (0.5,-0.5) {$K(\times)$};
\draw[black, very thick, ->] (3.4,0.6) -- (3,1);
\draw[black, very thick] (3.6,0.4) -- (4,0);
\draw[black, very thick, ->] (3,0) -- (4,1);
\node at (3.5,-0.5) {$K(\times^+)$};
\draw[black, very thick] (6,0) -- (6.4,0.4);
\draw[black, very thick, ->] (7,0) -- (6,1);
\draw[black, very thick, ->] (6.6,0.6) -- (7,1);
\node at (6.5,-0.5) {$K(\times^-)$};
\end{tikzpicture}
\end{center} \caption{The Vassiliev resolution.} \label{fig8}
\end{figure}
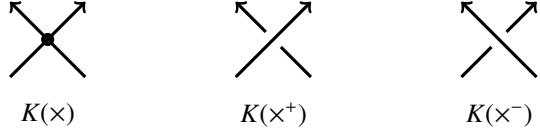

\begin{definition} {\rm
\emph{A singular crossing} is a crossing in a knotoid diagram decorated with a black dot, see \cite{MLK21}. An oriented \emph{singular planar knotoid diagram} is an oriented planar knotoid diagram that has finitely many singular crossings, as shown in Figure \ref{fig8}, see also~\cite{Bat23}. }
\end{definition}

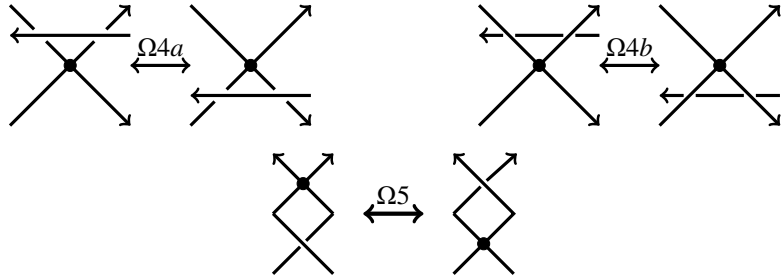
\begin{figure}[htbp]
	\centering
	\vspace{8pt}
	\begin{tikzpicture}[scale=0.8]
		\draw[black, very thick] (0,0)--(1.4,1.4);
		\draw[black, very thick,->] (1.6,1.6)--(2,2);
		\draw[black, very thick,<-] (0,1.5)--(2,1.5);
		\draw[black, very thick] (0,2)--(0.4,1.6);
		\draw[black, very thick,->](0.6,1.4)--(2,0);
		\draw[black, very thick,<->] (2,1)--(3,1);
		\node[above] at (2.5,1) {$\Omega4a$};
		\draw[black, very thick,<-] (3,0.5)--(5,0.5);
		\draw[black, very thick] (3,0)--(3.4,0.4);
		\draw[black, very thick,->] (3.6,0.6)--(5,2);
		\draw[black, very thick] (3,2)--(4.4,0.6);
		\draw[black, very thick,->](4.6,0.4)--(5,0);
		\filldraw[black] (1,1) circle (.1);
		\filldraw[black] (4,1) circle (.1);
	\end{tikzpicture}
	\qquad \qquad \qquad
	\begin{tikzpicture}[scale=0.8]
		\draw [black, very thick,->](0,1)--(2,3);
		\draw[black, very thick,->](0,3)--(2,1);
		\filldraw[black] (1,2) circle (.1);
		\filldraw[black] (4,2) circle (.1);
		\draw[black, very thick,<-](0,2.5)--(0.4,2.5);
		\draw[black, very thick](0.6,2.5)--(1.4,2.5);
		\draw[black, very thick](1.6,2.5)--(2,2.5);
		\draw[black, very thick,<->](2,2)--(3,2);
		\node[above] at (2.5,2) {$\Omega4b$};
		\draw[black, very thick,->](3,1)--(5,3);
		\draw[black, very thick,<-](3,1.5)--(3.4,1.5);
		\draw[black, very thick](3.6,1.5)--(4.4,1.5);
		\draw[black, very thick](4.6,1.5)--(5,1.5);
		\draw[black, very thick,->](3,3)--(5,1);
	\end{tikzpicture}
	
\bigskip
		\begin{tikzpicture}[scale=0.8]
		\draw[black, very thick,->](0,2)--(1,3);
		\draw[black, very thick](0,2)--(1,1);
		\draw[black, very thick,<-](0,3)--(1,2);
		\filldraw[black] (0.5,2.5) circle (.1);
		\draw[black, very thick](0,1)--(0.45,1.45);
		\draw[black, very thick](0.55,1.55)--(1,2);
		\draw[black, ultra thick,<->](1.5,2)--(2.5,2);
		\node[above] at (2,2){$\Omega5$};
		\draw[black, very thick,->](3.,1)--(4.,2)--(3,3);
		\draw[black, very thick,->](3.55,2.55)--(4.05,3);
		\draw[black, very thick] (3,2)--(3.45,2.45);
		\draw[black, very thick](3,2)--(4,1);
		\filldraw[black] (3.5,1.5) circle (.1);
	\end{tikzpicture}
	\caption{Singular Reidemeister moves $\Omega4a$, $\Omega4b$ and $\Omega5$.}
		\label{fig9}
\end{figure}

Two oriented singular planar knotoid diagrams are considered \emph{equivalent }if one can be obtained from the other by a finite sequence of oriented Reidemeister moves, singular Reidemeister moves as shown in Figure~\ref{fig9}, together with planar isotopies. An oriented singular planar knotoid is an equivalence class of oriented singular planar knotoid diagrams. The forbidden moves $\Phi_+$ and $\Phi_-$ of planar knotoids are also forbidden for singular planar knotoids.

Let $v$ be an invariant of a planar  knotoid, and let $v$ take values in an abelian group $G$. If $\times$ is a singular crossing of a planar knotoid $K(\times)$, then $K(\times^+)$ and $K(\times^-)$ correspond to the case where the singular crossing is a positive crossing and a negative crossing, respectively as shown in Figure \ref{fig8}. A knotoid invariant $v$ can be extended to a singular knotoid invariant by using the Vassiliev skein relation:
\begin{equation}
v(K(\times))=v(K(\times^+))-v(K(\times^-)). \label{eqn:Vas}
\end{equation}

If $v$ satisfies the condition that $v(K)\ne 0$ for a singular knotoid $K$ with $n$ singular points and $v$ vanishes for singular knotoids with more than $n$ singular points, then $v$ is called a Vassiliev invariant of \emph{order} $n$, see \cite{MLK21}.

\subsection{Crossing change and Gordian distance}
\begin{figure}[h]
\begin{center}
\begin{tikzpicture}
\draw[black, very thick] (0.4,0.6) -- (0,1);
\draw[black, very thick] (0.6,0.4) -- (1,0);
\draw[black, very thick] (0,0) -- (1,1);
\draw[black,very thick, <->] (2.5,0.5) -- (3.5,0.5);
\draw[black, very thick] (5,0) -- (5.4,0.4);
\draw[black, very thick] (6,0) -- (5,1);
\draw[black, very thick] (5.6,0.6) -- (6,1);
\end{tikzpicture}
\end{center} \caption{A crossing change.} \label{fig10}
\end{figure}
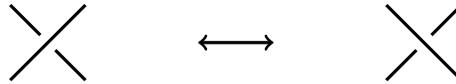

For a knot diagram, an unknotting operation is an operation that exchanges the over- and the undercrossing at a double point of the diagram. There are several unknotting operations on knot diagrams, and one of the simplest is a crossing change shown in Figure~\ref{fig10}. The unknotting number of a knot is the minimum number of unknotting operations required to reach a trivial knot among all diagrams of the knot, see~\cite{Mu85}.
Similar to classical knots, there is an unknotting operation on knotoid diagrams as well, which is called crossing change, see Figure~\ref{fig10}. It is well-known that any knot diagram can be transformed into a trivial knot diagram through crossing changes. But the analogue statement for knotoid diagrams may not hold, that is some knotoid diagrams can not be transformed into a trivial planar knotoid diagram through crossing changes. Two planar knotoid diagrams are said to be \emph{homotopic}, if they can be related by a finite sequence of crossing changes and Reidemeister moves.

\begin{remark} {\rm
There are more than one homotopy classes of planar knotoids. Indeed, by calculating the polynomial $X_p(K^R)$ defined in~\cite{Bat22}, of $K_1$ and $K_2$, presented in Figure~\ref{fig11}, one can see that both $K_1$ and $K_2$ are not homotopic to the trivial knotoid.}
\end{remark}
	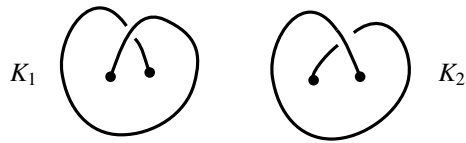
\begin{figure}[htbp]
	\begin{center}
\tikzset{every picture/.style={line width=1pt}}
\begin{tikzpicture}[x=0.6pt,y=0.6pt,yscale=-1,xscale=1]
	\draw[black, very thick]   (134.67,601.34) .. controls (139.33,589.34) and (148,556) .. (168.67,565.34) .. controls (189.33,574.67) and (192.67,586.67) .. (187.33,606) .. controls (182,625.34) and (161.33,638) .. (141.33,638) .. controls (121.33,638) and (107.33,622) .. (104,598.67) .. controls (100.67,575.34) and (124,539.34) .. (144.67,569.34) ;
	\filldraw[black] (134.67,601.34)  circle (1.5pt);
	\draw[black, very thick]   (150,577.34) .. controls (155.33,582.67) and (156,590.67) .. (159.33,598.67) ;
	\filldraw[black] (159.33,598.67) circle (1.5pt);
	\draw[black, very thick] (288,573.56) .. controls (311.33,554.89) and (327.33,585.56) .. (322,604.89) .. controls (316.67,624.23) and (294,640.67) .. (274,640.67) .. controls (254,640.67) and (240,624.67) .. (236.67,601.34) .. controls (233.33,578) and (256.67,542) .. (277.33,572) ;
	\draw[black, very thick]   (277.33,572) .. controls (282.67,577.34) and (288.67,593.34) .. (292,601.34) ;
	\filldraw[black] (292,601.34) circle (1.5pt);
	\draw[black, very thick]   (278,581.45) .. controls (270,587.56) and (259.33,600.89) .. (262.67,603.56) ;
	\filldraw[black] (262.67,603.56) circle (1.5pt);
	\node at (80,600) {$K_1$};
	\node at (350,600) {$K_2$};	
\end{tikzpicture}
	\caption{Knotoids $K_1$ and $K_2$.}  \label{fig11}
	\end{center}
\end{figure}	

\begin{definition} {\rm
Let $K$ and $K'$ be two homotopic planar knotoids, then the \emph{Gordian distance} denoted by $d_G(K,K')$ is the minimum number of crossing changes required to transform a diagram of $K$ into a diagram of $K'$. If $K$ and the trivial planar knotoid $\bf 0$ are homotopic, then $d_G(K,\bf 0)$ is the \emph{unknotting number} of $K$. }
\end{definition}

The introduction of a Gordian distance for an unknotting operation  allows to associate a simplicial complex with a set of classical or virtual knots. Namely, the Gordian complex $\mathcal G$ of knots is a simplicial complex defined by as follows. The vertex set of $\mathcal G$ consists of all the isotopy classes of oriented knots in $S^3$ and a family of $n+1$ vertices $\{ K_0, K_1, \ldots, K_n\}$ spans an $n$-simplex if and only if the Gordian distance $d_G(K_i, K_j) = 1$ for $0 \leq i, j \leq n$, $i \neq j$. For example, in~\cite{GPV}, Gill, Prabhakar and Vesnin showed the existence of an arbitrarily high-dimensional simplex in both the Gordian complex by regional crossing change for classical knots and the Gordian complex by arc shift move for virtual knots.
 In~\cite{KGPV}, Kaur, Gill, Prabhakar and Vesnin studied the Gordian distance between welded knots by twist move and properties of the corresponding Gordian complex. The natural question arises to study Gordian complexes corresponding to different distances between knotoids.

\section{Proofs of theorems} \label{sec3}

\subsection{Proof of Theorem~\ref{thm-inv}} \label{subsec3.1}

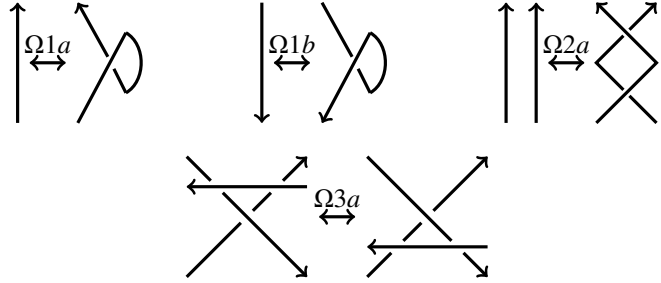
\begin{figure}[htbp]
	\centering
	\vspace{8pt}
	\begin{tikzpicture}[scale=0.8]
		\draw[black, very thick,->](0,1)--(0,3);
		\draw[black, very thick,<->](0.2,2)--(0.8,2);
		\node[above] at (0.5,2){$\Omega1a$};
		\draw[black, very thick, <-] (1,3)--(1.5,2.1);
		\draw[black, very thick] (1.6,1.9)--(1.8,1.5);
		\draw[black, very thick] (1,1)--(1.8,2.5);
		\draw[black, very thick] (1.8,1.5) to [out=30,in=330] (1.8,2.5);
		\node[above] at (0.5,0.5){};
	\end{tikzpicture}
	\qquad
	\qquad
	\begin{tikzpicture}[scale=0.8]
		\draw[black, very thick,<-](0,1)--(0,3);
		\draw[black, very thick,<->](0.2,2)--(0.8,2);
		\node[above] at (0.5,2){$\Omega1b$};
		\draw[black, very thick] (1,3)--(1.5,2.1);
		\draw[black, very thick] (1.6,1.9)--(1.8,1.5);
		\draw[black, very thick, <-] (1,1)--(1.8,2.5);
		\draw[black, very thick] (1.8,1.5) to [out=30,in=330] (1.8,2.5);
		\node[above] at (0.5,0.5){};
	\end{tikzpicture}
	\qquad
	\qquad
	\begin{tikzpicture}[scale=0.8]
		\draw[black, very thick,->](0,1)--(0,3);
		\draw[black, very thick,->](0.5,1)--(0.5,3);
		\draw[black, very thick,<->](0.7,2)--(1.3,2);
		\node[above] at (1,2){$\Omega2a$};
		\draw[black, very thick,->](1.5,1)--(2.5,2)--(1.5,3);
		\draw[black, very thick,->](2.05,2.55)--(2.5,3);
		\draw[black, very thick](1.95,1.55)--(1.5,2)--(1.95,2.45);
		\draw[black, very thick](2.5,1)--(2.05,1.45);
		\node[above]at(1,0.5){};
	\end{tikzpicture}
	\qquad
	\qquad
	\begin{tikzpicture}[scale=0.8]
		\draw[black, very thick](0,0)--(0.9,0.9);
		\draw[black, very thick](1.1,1.1)--(1.4,1.4);
		\draw[black, very thick,->](1.6,1.6)--(2,2);
		\draw[black, very thick,<-](0,1.5)--(2,1.5);
		\draw[black, very thick](0,2)--(0.4,1.6);
		\draw[black, very thick,->](0.6,1.4)--(2,0);
		\draw[black, very thick,<->](2.2,1)--(2.8,1);
		\node[above] at (2.5,1){$\Omega3a$};
		\draw[black, very thick,<-](3,0.5)--(5,0.5);
		\draw[black, very thick](3,0)--(3.4,0.4);
		\draw[black, very thick](3.6,0.6)--(3.9,0.9);
		\draw[black, very thick,->](4.1,1.1)--(5,2);
		\draw[black, very thick](3,2)--(4.4,0.6);
		\draw[black, very thick,->](4.6,0.4)--(5,0);
	\end{tikzpicture}
		\caption{Generating set of oriented Reidemeister moves.}
		\label{fig12}
\end{figure}

In~\cite{P10}, Polyak proved  that all oriented versions of Reidemeister moves for knot and link diagrams can be generated by a set of only four oriented Reidemeister moves, and that no less than four oriented Reidemeister moves generate them all. Thus, it is  suffices to prove that $H_D (t, y, z)$ is invariant under the oriented Reidemeister moves $\Omega 1a$, $\Omega 1b$, $\Omega 2a$ and $\Omega 3a$ of planar knotoid diagrams shown in Figure~\ref{fig12}.

Below we will consider the corresponding moves for Gauss diagrams.

Suppose $G(D)$ and $G(D')$ are Gauss diagrams of $D$ and $D'$, and $e$ and $e'$ are the crossings corresponding to $D$ and $D'$, respectively. Suppose the number of crossings in $D'$ is greater than or equal to the number of crossings corresponding to $D$.

\begin{lemma}\label{lem3.1}
$H_D(t,y,z)$ remains invariant under $\Omega1a$ and $\Omega1b$ moves.
\end{lemma}

\begin{proof}
Suppose that $D$ and $D'$ are two oriented knotoid diagrams which are related by a $\Omega1a$ move or a $\Omega1b$ move. Suppose $G(D)$ and $G(D')$ are the Gauss diagrams corresponding to $D$ and $D'$, respectively, and $c'$ is the new chord generated by the Reidemeister move. The Gauss diagrams corresponding to a $\Omega1a$ move and a  $\Omega1b$ move are shown in Figure~\ref{fig13}. 
	\begin{figure}[htbp] 
		\centering
		\begin{tikzpicture}[scale=0.7]
			\draw[black, thick](0,3.46) arc (120:240:2);
			\draw[black, thick,dashed] (2.73,2.73) arc (390:480:2);
			\draw[black, thick,dashed] (2.73,0.73) arc (330:240:2);
			\filldraw[black] (2.73,2.73) circle (.1);
			\filldraw[black] (2.73,0.73) circle (.1);		
		\end{tikzpicture}
		\quad
			\begin{tikzpicture}[scale=0.7]
			\draw [very thick, black,<->](0,3) --(2,3);
			\node [above] at(1,3){$\Omega1a$};
			\node [above] at(1.5,1){\ };
		\end{tikzpicture}
		\quad
		\begin{tikzpicture}[scale=0.7]
			 \draw[black, thick](0,3.46) arc (120:240:2);
			\draw[black, thick,dashed] (2.73,2.73) arc (390:480:2);
			\draw[black, thick,dashed] (2.73,0.73) arc (330:240:2);
			\filldraw[black] (2.73,2.73) circle (.1);
			\filldraw[black] (2.73,0.73) circle (.1);
			\draw [black, very thick,->](-0.28,3.26) --(-0.28,0.2);
			\node [above] at(-0.28,3.26){$+$};
			\node[right] at (-0.28,1.8){$c'$};
		\end{tikzpicture}
		\\
		\begin{tikzpicture}[scale=0.7]
			\draw[black, thick](0,3.46) arc (120:240:2);
			\draw[black, thick, dashed] (2.73,2.73) arc (390:480:2);
			\draw[black, thick, dashed] (2.73,0.73) arc (330:240:2);
			\filldraw[black] (2.73,2.73) circle (.1);
			\filldraw[black] (2.73,0.73) circle (.1);			
		\end{tikzpicture}
		\quad
		\begin{tikzpicture}[scale=0.7]
			\draw [black, very thick,<->](0,3) --(2,3);
			\node [above] at(1,3){$\Omega1b$};
			\node [above] at(1.5,1){\ };
		\end{tikzpicture}
		\quad
		\begin{tikzpicture}[scale=0.7]
			\draw[black, thick](0,3.46) arc (120:240:2);
			\draw[black, thick, dashed] (2.73,2.73) arc (390:480:2);
			\draw[black, thick, dashed] (2.73,0.73) arc (330:240:2);
			\filldraw[black] (2.73,2.73) circle (.1);
			\filldraw[black] (2.73,0.73) circle (.1);
			\draw [black, very thick,<-](-0.28,3.26) --(-0.28,0.2);
			\node [below]at(-0.28,0.2){$+$};
			\node[right] at (-0.28,1.8){$c'$};
		\end{tikzpicture}
		\caption{$\Omega$1a move and $\Omega$1b move.}
		\label{fig13}
	\end{figure}
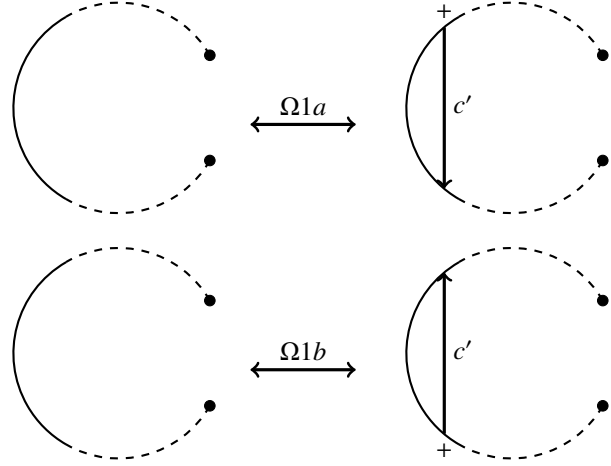
Since $c'$ is isolated and no other chords pass through it, then $d(c')=0$ and $\operatorname{Ind}_{c'}^n(z)=0$. Hence $H_{D'}(t,y,z) - H_D(t,y,z)=\operatorname{sgn}(c')(t^{ \operatorname{Ind}_{c'}^n(z)}-1)y^n=0$.
\end{proof}

\begin{lemma}\label{lem3.2}
$H_D(t,y,z)$ remains invariant under a $\Omega2a$ move.
\end{lemma}

\begin{proof}
Suppose that $D$ and $D'$ are two knotoid diagrams which are related by a $\Omega2a$ move. Suppose $G(D)$ and $G(D')$ are the Gauss diagrams corresponding to $D$ and $D'$, respectively, and $c_1,c_2$ are the new chords generated by the Reidemeister move. The Gauss diagrams corresponding to a $\Omega2a$ move are shown in Figure~\ref{fig14}. 
	\begin{figure}[htbp] 
	\centering
	\begin{tikzpicture}[scale=0.7]
		\draw[black, thick](2,3.46) arc (60:120:2);
		\draw[black, thick](0,0) arc (240:300:2);
		\draw[black, thick,dashed] (2.73,0.73) arc (330:300:2);
		\draw[black, thick,dashed] (2.73,2.73) arc (390:420:2);
		\draw[black, thick,dashed] (0,3.46) arc (120:240:2);
		\filldraw[black] (2.73,2.73) circle (.1);
		\filldraw[black] (2.73,0.73) circle (.1);
	\end{tikzpicture}
	\quad
	\begin{tikzpicture}[scale=0.7]
	\draw [black, very thick,<->](0,3) --(2,3);
	\node [above] at(1,3){$\Omega2a$};
	\node [above] at(1.5,1){\ };
    \end{tikzpicture}
	\quad
	\begin{tikzpicture}[scale=0.7]
		\draw[black, thick](2,3.46) arc (60:120:2);
		\draw[black, thick](0,0) arc (240:300:2);
		\draw[black, thick, dashed] (2.73,0.73) arc (330:300:2);
		\draw[black, thick, dashed] (2.73,2.73) arc (390:420:2);
		\draw[black, thick, dashed] (0,3.46) arc (120:240:2);
		\filldraw[black] (2.73,2.73) circle (.1);
		\filldraw[black] (2.73,0.73) circle (.1);
		\draw [black, very thick,->](1.68,3.6) --(0.32,-0.18);
		\draw [black, very thick,->](0.32,3.6) --(1.68,-0.18);
		\node [above] at(1.68,3.6){$+$};
		\node [above]at(0.32,3.6){$-$};
		\node[right] at (1.33,2.28){$c_1$};
		\node[left] at (0.83,2.28){$c_2$};
	\end{tikzpicture}
		\caption{$\Omega2a$ move.}
		\label{fig14}
\end{figure}
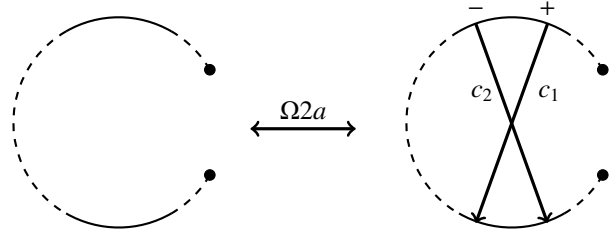

It is easy to see that $\operatorname{sgn}(c_1) = -\operatorname{sgn}(c_2)$ and $d(c_1)=d(c_2)$. For each $e\in C(G(D)) \setminus \{ c_1, c_2 \}$, denote it by $e'$ the chord in $C(G(D'))$ corresponding to $e$ after the $\Omega2a$ move. Since $c_1$ and $c_2$ pass through $e'$ in the same direction, while their signs are opposite, the affect of $c_1$ and $c_2$ offsets, which implies $d(e)=d(e')$.
			
Since $d(c_1)=d(c_2)$, then for any $n \in \mathbb N$ we have  $c_1\in C_n(G(D'), e')$ if and only if $c_2\in C_n(G(D'), e')$. Hence the affect of $c_1$ and $c_2$ offsets, which implies $\operatorname{Ind}_e^n(z) = \operatorname{Ind}_{e'}^n(z)$.
			
For chords $c_1$ and $c_2$, if $n = d(c_1) = d(c_2) = \gcd (d(c_1),d(c_2))$, then $\ell^n(c_1)\backslash \ell^n(c_2) = \{ c_2 \}$ and $r^n(c_2)\backslash r^n(c_1)=\{ c_1\}$. As a consequence, we get
$$
\operatorname{Ind}_{c_1}^n(z) - \operatorname{Ind}_{c_2}^n(z) = -\operatorname{sgn}(c_2) \, z^{ \phi _{c_1}(-d(c_2))} - \operatorname{sgn}(c_1) \, z^{ \phi _{c_2}(d(c_1))},
$$
where $\phi _{c_1}:\mathbb{Z}\rightarrow \mathbb{Z}_{|d(c_1)|}=\mathbb{Z}_n$ and $\phi _{c_2}:\mathbb{Z}\rightarrow \mathbb{Z}_{|d(c_2)|}=\mathbb{Z}_n$. Hence maps $\phi_{c_1}$ and $\phi_{c_2}$ coincide, $\phi _{c_1} = \phi _{c_2}$, and so $\operatorname{Ind}_{c_1}^n(z) = \operatorname{Ind}_{c_2}^n(z)$.

If $n\in\mathbb{N}$ and $n\ne d(c_1)$ (and $n \neq d(c_2)$ since $d(c_1) = d(c_2)$), then $\ell^n(c_1)=\ell^n(c_2)$ and $r^n(c_1)=r^n(c_2)$, which implies $\operatorname{Ind}_{c_1}^n(z) =  \operatorname{Ind}_{c_2}^n(z)$. Therefore,
$$
H_{D'}(t,y,z) - H_D(t,y,z) = \operatorname{sgn}(c_1) (t^{\operatorname{Ind}_{c_1}^n(z)}-1)y^n + \operatorname{sgn}(c_2) (t^{\operatorname{Ind}_{c_2}^n(z)}-1)y^n=0,
$$
which completes the proof of the Lemma.
\end{proof}

\begin{lemma}\label{lem3.3}
$H_D(t,y,z)$ remains invariant under a $\Omega3a$ move.
\end{lemma}	

\begin{proof}
For the oriented Reidemeister move $\Omega3a$, there are two types of Gauss diagrams corresponding to it: $\Omega3a$ and $\Omega3a'$. Suppose that $D$ and $D'$ are two knotoid diagrams which are related by an $\Omega3a$ move or $\Omega3a'$ move. Suppose $G(D)$ and $G(D')$ are the Gauss diagrams corresponding to $D$ and $D'$, respectively, and $c_1$, $c_2$, and $c_3$ in $G(D)$ corresponds one-to-one with $c'_1$, $c'_2$, and $c'_3$ in $G(D')$, see Figure~\ref{fig15}.
 \begin{figure}[htbp] 
\centering
\begin{tikzpicture}[scale=0.7]
\draw[black, thick](-0.1,0.05) arc (30:60:2);
\draw[black, thick](-2.73,0.73) arc (120:150:2);
\draw[black, thick](-2.73,-2.73) arc (240:300:2);
\draw[red, very thick,dashed] (-0.73,-2.73) arc (300:330:2);
\draw[red,very thick,dashed] (0.15,-0.25) arc (30:42:2);
\draw[green,very thick,dashed] (-0.73,0.73) arc (60:120:2);
\draw[brown,very thick,dashed] (-3.46,0) arc (150:240:2);
\node at (-4,-1.5) {$c$};
\node at (-1.8,1.5) {$b$};
\node at (0.3,0) {$a$};
\node at (0.1,-2.5) {$a$};
\filldraw[black] (0.15,-0.3) circle (.1);
\filldraw[black] (0,-2) circle (.1);
\draw [black, very thick,->](-2.73,-2.73) --(-0.73,0.73);
\draw [black, very thick,->](-0.73,-2.73) --(-2.73,0.73);
\draw [black, very thick,->](-3.3,0.2) --(-0.15,0.2);
\node [left] at(-3.46,0){$+$};
\node [below] at(-2.73,-2.73){$-$};
\node [below] at(-0.73,-2.73){$+$};
\node [below] at(-1.8,0.7){$c_3$};
 \node [left] at(-2,-0.5){$c_1$};
 \node [right] at(-1.3,-0.5){$c_2$};
 \end{tikzpicture}
	\quad
	\begin{tikzpicture}[scale=0.7]
		\draw [black, very thick,<->](0,3) --(2,3);
		\node [above] at(1,3){$\Omega3a$};
		\node [above] at(1.5,1){\ };
	\end{tikzpicture}
	\quad
	\begin{tikzpicture}[scale=0.7]
	\node at (-4,-1.5) {$c$};
	\node at (-1.8,1.5) {$b$};
         \node at (0.3,0) {$a$};
	\node at (0.1,-2.5) {$a$};
	\draw[black, thick](-0.1,0.05) arc (30:60:2);
	\draw[black, thick](-2.73,0.73) arc (120:150:2);
	\draw[black, thick](-2.73,-2.73) arc (240:300:2);
	\draw[red,very thick,dashed] (0.15,-0.25) arc (30:42:2);
	\draw[green, very thick, dashed] (-0.73,0.73) arc (60:120:2);
	\draw[purple, very thick, dashed] (-3.46,0) arc (150:240:2);
	\draw[red, very thick,dashed] (-0.73,-2.73) arc (300:330:2);
	\draw[purple, very thick,dashed] (0.15,-0.3) arc (30:40:2);
	\filldraw[black] (0.15,-0.3) circle (.1);
	\filldraw[black] (0,-2) circle (.1);
	\draw [black, very thick,->](-2.73,0.73) --(-0.73,0.73);
	\draw [black, very thick,->](-2.73,-2.73) --(-3.3,0.2);
	\draw [black, very thick,->](-0.73,-2.73) --(-0.15,0.2);
	\node [above] at(-3,0.5){$+$};
	\node [below] at(-2.73,-2.73){$+$};
	\node [below] at(-0.73,-2.73){$-$};
	\node [below] at(-1.8,0.6){$c'_3$};
	\node [left] at(-2.5,-0.5){$c'_1$};
	\node [right] at(-0.8,-0.5){$c'_2$};
\end{tikzpicture}
\\
\begin{tikzpicture}[scale=0.7]
	\node at (-4,-1.5) {$c$};
	\node at (-1.8,1.5) {$b$};
         \node at (0.3,0) {$a$};
	\node at (0.1,-2.5) {$a$};
	\draw[black, thick](-0.1,0.05) arc (30:60:2);
	\draw[black, thick](-2.73,0.73) arc (120:150:2);
	\draw[black, thick](-2.73,-2.73) arc (240:300:2);
			\draw[red, very thick,dashed] (-0.73,-2.73) arc (300:330:2);
		\draw[red,very thick,dashed] (0.15,-0.25) arc (30:42:2);
		\draw[green,very thick,dashed] (-0.73,0.73) arc (60:120:2);
		\draw[brown,very thick,dashed] (-3.46,0) arc (150:240:2);
	\filldraw[black] (0.15,-0.3) circle (.1);
	\filldraw[black] (0,-2) circle (.1);
	\draw [black, very thick,->](-2.73,0.73) --(-0.73,0.73);
	\draw [black, very thick,<-](-2.73,-2.73) --(-3.3,0.2);
	\draw [black, very thick,->](-0.73,-2.73) --(-0.15,0.2);
	\node [above] at(-3.,0.5){$-$};
	\node [left] at(-3.4,0.2){$+$};
	\node [below] at(-0.73,-2.73){$+$};
	\node [below] at(-1.8,0.6){$c_3$};
	\node [left] at(-2.5,-0.5){$c_1$};
	\node [right] at(-0.8,-0.5){$c_2$};
\end{tikzpicture}
\quad
\begin{tikzpicture}[scale=0.7]
	\draw [black, very thick,<->](0,3) --(2,3);
	\node [above] at(1,3){$\Omega3a'$};
	\node [above] at(1.5,1){\ };
\end{tikzpicture}
\quad
	\begin{tikzpicture}[scale=0.7]
		\node at (-4,-1.5) {$c$};
	\node at (-1.8,1.5) {$b$};
         \node at (0.3,0) {$a$};
	\node at (0.1,-2.5) {$a$};
	\draw[black, thick](-0.1,0.05) arc (30:60:2);
	\draw[black, thick](-2.73,0.73) arc (120:150:2);
	\draw[black, thick](-2.73,-2.73) arc (240:300:2);
			\draw[red, very thick,dashed] (-0.73,-2.73) arc (300:330:2);
		\draw[red,very thick,dashed] (0.15,-0.25) arc (30:42:2);
		\draw[green,very thick,dashed] (-0.73,0.73) arc (60:120:2);
		\draw[brown,very thick,dashed] (-3.46,0) arc (150:240:2);
	\filldraw[black] (0.15,-0.3) circle (.1);
	\filldraw[black] (0,-2) circle (.1);
	\draw [black, very thick,->](-2.73,-2.73) --(-0.73,0.73);
	\draw [black, very thick,<-](-0.73,-2.73) --(-2.73,0.73);
	        \draw [black, very thick,->](-3.3,0.2) --(-0.15,0.2);
	\node [above] at(-2.73,0.73){$+$};
	\node [left] at(-3.46,0.2){$-$};
	\node [below] at(-2.73,-2.73){$+$};
	\node [below] at(-1.8,0.9){$c'_3$};
	\node [left] at(-2,-0.5){$c'_1$};
	\node [right] at(-1.3,-0.5){$c'_2$};
\end{tikzpicture}
\caption{$\Omega3a$ move and $\Omega3a'$ move.} \label{fig15}
\end{figure}
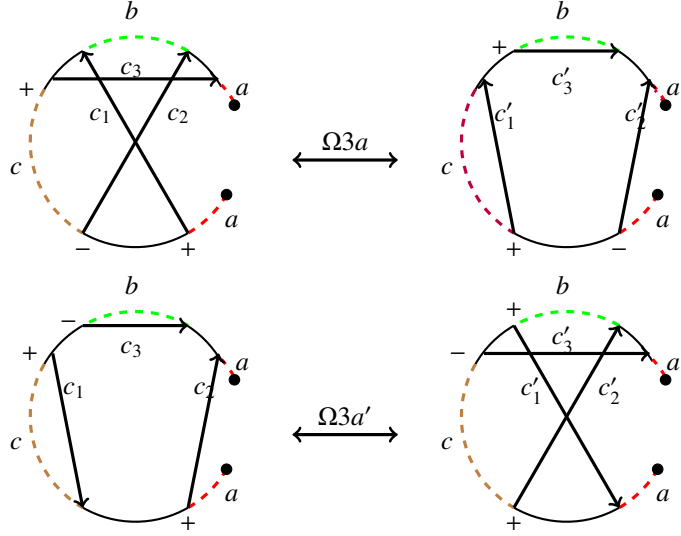

Let us introduce some notations consisting which of two letters and a sign as follows. First divide the Gauss diagram $G(D)$ into three regions $a$, $b$ and $c$ as shown in the Figure~\ref{fig15}. Mark the $a$ region in red,  the $b$ region in green, and the $c$ region in brown. Then write out the regions where the two endpoints of a chord are located. The first letter is the region where the tail of the chord is located, and the second letter is the region where the head of the chord is located. Finally add the sign of the chord. These notations indicate the number of chords whose tails are in the region of the first letter, whose heads are in the region of the second letter and whose signs are the corresponding ones in the notation. For example, $ab^+$ denotes the number of chords with tails in $a$, heads in $b$, and positive signs. Then we can see that
$$	
\begin{gathered}
d(c_1)=cb^+ + ca^+ - cb^- - ca^- - bc^+ - ac^+ + bc^- + ac^- +1-1, \cr
d(c_2)=ba^+ + ca^+ - ba^- - ca^- - ab^+ - ac^+ + ab^- + ac^- + 1-1, \cr
d(c_3)=bc^+ + ba^+ - bc^- - ba^- - ab^+ - cb^+ + ab^- + cb^- + 1-1.
\end{gathered}
$$	
It is easy to see that $d(c_1)+d(c_3)=d(c_2)$, whence
$$
\gcd(d(c_1),d(c_2)) = \gcd(d(c_2),d(c_3)) = \gcd(d(c_1),d(c_3)).
$$
It is clear that
$$
\operatorname{sgn}(c_1) = \operatorname{sgn}(c'_1) =+1, \quad \operatorname{sgn}(c_2) = \operatorname{sgn}(c'_2)=-1, \quad \operatorname{sgn}(c_3) = \operatorname{sgn}(c'_3)=+1,
$$
hence $d(c_i)=d(c'_i)$ for $i=1,2,3$. If $e \in C(G(D)) \setminus \{ c_1, c_2, c_3 \}$ then  $d(e)=d(e')$.
	
For any natural number $n$, due to the unchanged direction and sign of the chord before and after the move, we see that $\operatorname{Ind}_e^n(z) = \operatorname{Ind}_{e'}^n(z)$.
	
Consider chords $c_1,c_2,c_3$ and $c'_1,c'_2,c'_3$. If
$$
n=\gcd (d(c_1),d(c_2)) = \gcd (d(c_2),d(c_3)) = \gcd (d(c_1),d(c_3)),
$$
then
\begin{alignat}{4}
\operatorname{Ind}_{c_1}^n(z) - \operatorname{Ind}_{c'_1}^n(z)
&= \operatorname{sgn}(c_2) \, z^{\phi _{c_1}(d(c_2))} + \operatorname{sgn}(c_3) \, z^{\phi _{c_1}(d(c_3))}
=-z^{\phi_{c_1}(d(c_2))} + z^{\phi_{c_1}(d(c_2)-d(c_1))}\nonumber \\
&=- z^{\phi_{c_1}(d(c_2))} + z^{\phi_{c_1}(d(c_2))} = 0,\nonumber
\end{alignat}
where $\phi _{c_1}:\mathbb{Z}\rightarrow \mathbb{Z}_{|d(c_1)|}$. Analogously,
\begin{alignat}{4}
	\operatorname{Ind}_{c_2}^n(z) - \operatorname{Ind}_{c'_2}^n(z)
	&=\operatorname{sgn}(c_3) \, z^{\phi _{c_2}(d(c_3))} - \operatorname{sgn}(c_1) \, z^{\phi _{c_2}(-d(c_1))}
	= z^{\phi_{c_2}(d(c_2)-d(c_1))} - z^{\phi_{c_2}(-d(c_1))}\nonumber\\
	&= z^{\phi_{c_2}(-d(c_1))} - z^{\phi_{c_2}(-d(c_1))}
	=0\nonumber,
\end{alignat}
where $\phi _{c_2}:\mathbb{Z}\rightarrow \mathbb{Z}_{|d(c_2)|}$. And
\begin{alignat}{4}
	\operatorname{Ind}_{c_3}^n(z) - \operatorname{Ind}_{c'_3}^n(z)
	&=-\operatorname{sgn}(c_2) z^{\phi _{c_3}(-d(c_2))} - \operatorname{sgn}(c_1) z^{\phi _{c_3}(-d(c_1))}
	=z^{\phi_{c_3}(-d(c_3)-d(c_1))} - z^{\phi_{c_3}(-d(c_1))} \nonumber\\
	&= z^{\phi_{c_3}(-d(c_1))} - z^{\phi_{c_3}(-d(c_1))} = 0 \nonumber,
\end{alignat}
where $\phi _{c_3}:\mathbb{Z}\rightarrow \mathbb{Z}_{|d(c_3)|}$.

	If $n \ne \gcd (d(c_1),d(c_2))$, then it is easy to see that $\operatorname{Ind}_{c_i}^n(z) = \operatorname{Ind}_{c'_i}^n(z)$, where $i=1,2,3$.
	
	Summarizing, for any $n \in \mathbb N$ we get $\operatorname{Ind}_{c_i}^n(z) = \operatorname{Ind}_{c'_i}^n(z)$, where $i=1,2,3$. Therefore,
$$
H_{D'}(t,y,z) - H_D(t,y,z)=\sum_{i=1}^{3} \operatorname{sgn}(c_i) (t^{\operatorname{Ind}_{c_i}^n(z)}-1)y^n-\sum_{i=1}^{3} \operatorname{sgn}(c'_i)(t^{\operatorname{Ind}_{c'_i}^n(z)}-1)y^n=0.
$$
Whence $H_{D'}(t,y,z) = H_D(t,y,z)$.

The proof of invariance of $H_D(t,y,z)$ under an $\Omega3a'$ move is similar.
\end{proof}

Summarizing, Theorem~\ref{thm-inv} holds from Lemmas \ref{lem3.1}-\ref{lem3.3}.

\subsection{Proof of Theorem~\ref{thm-prop}}

\underline{Proof of point (1)}.
	The Gauss diagrams corresponding to $D$ and $-D$ are shown in Figure~\ref{fig16}. Let $c$ be a chord in $D$ corresponding to $c'$ in $-D$. $G(D)$ and $G(-D)$ are the Gauss diagrams corresponding to $D$ and $-D$, respectively. It can be seen that $G(-D)$ is the Gauss diagram obtained by the $\pi$-rotation of $G(D)$ around the horizontal axis. Since $\operatorname{sgn}(c )= \operatorname{sgn}(c')$ and $d(c)=-d(c')$, as well as $r^n(c) = \ell^n(c')$ and $\ell^n(c)=r^n(c')$, we get  $\phi_{c}=\phi_{c'}$ and $\operatorname{Ind}_c^n(z)=-\operatorname{Ind}_{c'}^n(z)$. Therefore, $H_D(t, y, z)=H_{-D}(t^{-1}, y, z)$.
	
 \begin{figure}[htbp] 
	\centering
	\begin{tikzpicture}[scale=0.7]
		\draw[black, thick](0,0) arc (30:330:2);
		\filldraw[black] (0,0) circle (.1);
		\filldraw[black] (0,-2) circle (.1);
		\draw [black, very thick,->](-0.21,0.28) -- (-3.25,0.28);
		\node [below]at(-2.5,0.28){$c_2$};
		\draw [black, very thick,->](-3.25,-2.28) -- (-0.21,-2.28);
		\node [above]at(-2.5,-2.28){$c_1$};
		\draw [black, very thick,->](-1.7,1) -- (-1.7,-3);
		\node [above] at(-1.5,1){$\operatorname{sgn}(c)$};
		\node [below]at(-1.5,-3){\,};
		\node[right] at (-1.5,-1){$c$};
		\node [below]at(-1.5,-4){$G(D)$};
	\end{tikzpicture}
	\quad\quad
		\begin{tikzpicture}[scale=0.7]
		\draw[black, thick](0,0) arc (30:330:2);
		\filldraw[black] (0,0) circle (.1);
		\filldraw[black] (0,-2) circle (.1);
		\draw [black, very thick,<-](-0.21,0.28) -- (-3.25,0.28);
		\node [below]at(-2.5,0.28){$c'_1$};
		\draw [black, very thick,<-](-3.25,-2.28) -- (-0.21,-2.28);
		\node [above]at(-2.5,-2.28){$c'_2$};
		\draw [black, very thick,<-](-1.7,1) -- (-1.7,-3);
		\node [above] at(-1.5,1){\,};
		\node [below]at(-1.5,-3){$\operatorname{sgn}(c)$};
		\node[right] at (-1.5,-1){$c'$};
		\node [below]at(-1.5,-4){$G(-D)$};
	\end{tikzpicture}
		\caption{The Gauss diagrams of $D$ and $-D$.}
		\label{fig16}
\end{figure}

\underline{Proof of point (2)}.
	The Gauss diagrams corresponding to $D$ and $D^*$ are shown in Figure~\ref{fig17}. Let $c$ be a chord in $D$ corresponding to $c'$ in $D^*$. $G(D)$ and $G(D^*)$ are the Gauss diagrams corresponding to $D$ and $D^*$, respectively. It can be seen that $G(D^*)$ is obtained by changing the sign and direction of each chord in $G(D)$. Hence $d(c)=-d(c')$, and $r^n(c)=r^n(c')$, $\ell^n(c)=\ell^n(c')$. But the signs of the chords in the sets are different, so $\phi_{c}=\phi_{c'}$ and $\operatorname{Ind}_c^n(z) = - \operatorname{Ind}_{c'}^n(z^{-1})$, hence $H_D(t,y,z)=-H_{D^*}(t^{-1},y,z^{-1})$. Therefore, $H_D(t,y,z)$ can distinguish whether the knotoid is chiral or not.
 \begin{figure}[htbp] 
	\centering
	\begin{tikzpicture}[scale=0.7]
		\draw[black, thick](0,0) arc (30:330:2);
		\filldraw[black] (0,0) circle (.1);
		\filldraw[black] (0,-2) circle (.1);
		\draw [black, very thick,->](-0.21,0.28) -- (-3.25,0.28);
		\node [below]at(-2.5,0.28){$c_1$};
		\draw [black, very thick,->](-3.25,-2.28) -- (-0.21,-2.28);
		\node [above]at(-2.5,-2.28){$c_2$};
		\draw [black, very thick,->](-1.7,1) -- (-1.7,-3);
		\node [above] at(-1.5,1){$\operatorname{sgn} (c)$};
		\node [below]at(-1.5,-3){\,};
		\node[right] at (-1.5,-1){$c$};
		\node [below]at(-1.5,-4){$G(D)$};
	\end{tikzpicture}
	\quad\quad
	\begin{tikzpicture}[scale=0.7]
		\draw[black, thick](0,0) arc (30:330:2);
		\filldraw[black] (0,0) circle (.1);
		\filldraw[black] (0,-2) circle (.1);
		\draw [black, very thick,<-](-0.21,0.28) -- (-3.25,0.28);
		\node [below]at(-2.5,0.28){$c'_1$};
		\draw [black, very thick,<-](-3.25,-2.28) -- (-0.21,-2.28);
		\node [above]at(-2.5,-2.28){$c'_2$};
		\draw [black, very thick,<-](-1.7,1) -- (-1.7,-3);
		\node [above] at(-1.5,1){\,};
		\node [below]at(-1.5,-3){$\operatorname{sgn}(c)$};
		\node[right] at (-1.5,-1){$c'$};
		\node [below]at(-1.5,-4){$G(D^*)$};
	\end{tikzpicture}
	\caption{The Gauss diagrams of $D$ and $D^*$.}
		\label{fig17}
\end{figure}

\underline{Proof of point (3)}.
Since $D$ is a zero height planar knotoid diagram, for any $c \in C(G(D))$ we have $d(c)=0$, so
$\operatorname{Ind}_c^n(z)=0$, whence $H_D(t,y,z)=0$.

Thus, Theorem~\ref{thm-prop} is proved. \hfill $\Box$

\subsection{Proof of Theorem~\ref{thm-Vas}}	

To prove that $H_K(t,y,z)$ is a Vassiliev invariant of order one for
planar knotoids, we will first show that $H_{D(\times_1,\times_2)}(t,y,z)=0$, where $D(\times_1,\times_2)$ is a singular planar knotoid diagram with two singular crossings $\times_1$ and $\times_2$. With the Vassiliev skein relation (\ref{eqn:Vas}),
\begin{alignat}{2}
	H_{D(\times_1,\times_2)}(t,y,z)
	&=H_{D(\times_1,\times_2^+)}(t,y,z) - H_{D(\times_1,\times_2^-)}(t,y, z)\nonumber\\
	&=H_{D(\times_1^+,\times_2^+)}(t,y,z)-H_{D(\times_1^+,\times_2^-)}(t,y,z)-H_{D(\times_1^-,\times_2^+)}(t,y,z)+H_{D(\times_1^-,\times_2^-)}(t,y,z)\nonumber.
\end{alignat}
Let $c$ be a classical crossing of $D(\times_1,\times_2)$ and $\varepsilon$ denotes "$+$" or "$-$". We see that
$$
d(c,D(\times_1^+,\times_2^+))=d(c,D(\times_1^+,\times_2^-))=d(c,D(\times_1^-,\times_2^+))=d(c,D(\times_1^-,\times_2^-))
$$
and for $n \in \mathbb N$ we get
$$
\operatorname{Ind}_{c,D(\times_1^\varepsilon,\times_1^\varepsilon)}^n(z) = \operatorname{Ind}_{c,D(\times_1^+,\times_1^+)}^n(z).
$$

For the two singular crossings $\times_1$ and $\times_2$ we get
$$
\begin{gathered}
d(\times_1^+,D(\times_1^+,\times_2^\varepsilon))=-d(\times_1^-,D(\times_1^+,\times^\varepsilon_2)), \qquad
d(\times_2^+,D(\times_1^\varepsilon,\times_2^+))=-d(\times_2^-,D(\times_1^\varepsilon,\times_2^-)), \cr
\operatorname{Ind}_{\times_1^+,D(\times_1^+,\times_2^+)}^n(z) = \operatorname{Ind}_{\times_1^+,D(\times_1^+,\times_2^-)}^n(z), \qquad
\operatorname{Ind}_{\times_1^-,D(\times_1^-,\times_2^+)}^n(z) = \operatorname{Ind}_{\times_1^-,D(\times_1^-,\times_2^-)}^n(z), \cr
\operatorname{Ind}_{\times_2^+,D(\times_1^+,\times_2^+)}^n(z) = \operatorname{Ind}_{\times_2^+,D(\times_1^-,\times_2^+)}^n(z),\qquad
\operatorname{Ind}_{\times_2^-,D(\times_1^+,\times_2^-)}^n(z) = \operatorname{Ind}_{\times_2^-,D(\times_1^-,\times_2^-)}^n(z).
\end{gathered}
$$

Then the coefficient of $y^n$ is
\begin{alignat}{2}
&\sum_{n\in \mathbb{N}}\left [ \left( t^{\operatorname{Ind}_{ \times_1^+,D(\times_1^+,\times_2^+)} ^n(z)}-1 \right) - \left ( t^{\operatorname{Ind}_{ \times_1^+,D(\times_1^+,\times_2^-) }^n(z)}-1 \right )\right ] \nonumber\\
&+\sum_{n\in \mathbb{N}}\left [ \left ( t^{\operatorname{Ind}_{ \times_1^-,D(\times_1^-,\times_2^+) } ^n(z)}-1 \right ) - \left ( t^{\operatorname{Ind}_{ \times_1^-,D(\times_1^-,\times_2^-) }^n(z)}-1 \right )\right ] \nonumber\\
&+\sum_{n\in \mathbb{N}}\left [ \left ( t^{\operatorname{Ind}_{ \times_2^+,D(\times_1^+,\times_2^+) } ^n(z)}-1 \right ) - \left (t^{\operatorname{Ind}_{ \times_2^+,D(\times_1^-,\times_2^+) }^n(z)}-1 \right )\right ] \nonumber\\
&+\sum_{n\in \mathbb{N}}\left [ \left ( t^{\operatorname{Ind}_{ \times_2^-,D(\times_1^+,\times_2^-) } ^n(z)}-1 \right ) - \left ( t^{\operatorname{Ind}_{\times_2^-,D(\times_1^-,\times_2^-)  }^n(z)}-1 \right) \right ] =0\nonumber.
\end{alignat}
Therefore, $H_{D(\times_1,\times_2)}(t,y,z)=0$.

In Example~\ref{exam4.3}, we calculate that
$$
H_{K(c)}(t,y,z)=(t^{-z}-1)y+(t^{-1}-1)y^2+(t^z-1)y+(t-1)y^2\ne 0,
$$
where $K(c)$ is a singular planar knotoid with one singular crossing $c$ as shown in Figure~\ref{fig22}. So $H_K(t,y,z)$ is a Vassiliev invariant of order one for planar knotoids. 

Thus, Theorem~\ref{thm-Vas} is proved. \hfill $\Box$

\subsection{Proof of Theorem~\ref{thm-dist}}

We start with the following lemma.

\begin{lemma}\label{lem5.1-CCmove}
Let $K$ and $K'$ be two planar knotoid diagram differ by a crossing change move at the crossing corresponding to a chord $c$. Then
$$
H_K(t,y,z)-H_{K'}(t,y,z) = \sum_{n\in\mathbb{N}} \varepsilon \left( t^{\operatorname{Ind}_c^n(z)} + t^{ -\operatorname{Ind}_c^n(z^{-1})}-2\right) y^n,
$$
where $\varepsilon = \operatorname{sgn} (c) \in \left\lbrace +1,-1\right\rbrace $.
\end{lemma}

\begin{proof}
The local Gauss diagrams of $K$ and $K'$ are shown in Figure \ref{fig18}.
	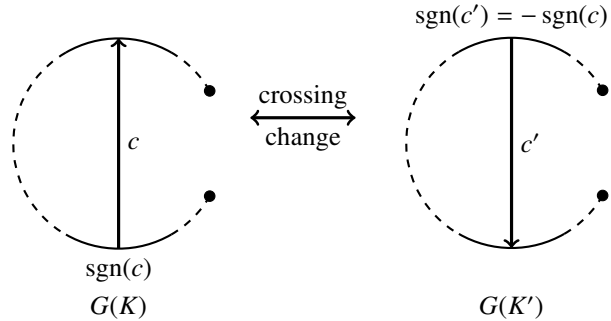
\begin{figure}[htbp] 
	\centering
	\begin{tikzpicture}[scale=0.7]
		\draw[black, thick](2,3.46) arc (60:120:2);
		\draw[black, thick](0,0) arc (240:300:2);
		\draw[thick,dashed] (2,0) arc (300:330:2);
		\draw[thick,dashed] (2.73,2.73) arc (390:420:2);
		\draw[black, thick,dashed] (0,3.46) arc (120:240:2);
		\filldraw[black] (2.73,2.73) circle (.1);
		\filldraw[black] (2.73,0.73) circle (.1);
		\draw [black, very thick,<-](1,3.73) --(1,-0.27);
		\node[right] at (1,1.73){$c$};
		\node[above] at (1,3.73){\,};
		\node[below] at (1,-0.27){$\operatorname{sgn}(c)$};
			\node[below] at (1,-1){$G(K)$};
	\end{tikzpicture}
	\quad
	\begin{tikzpicture}[scale=0.7]
		\draw [black, very thick,<->](0,3) --(2,3);
		\node [above] at(1,3){crossing};
		\node [below] at(1,3){change};
		\node [above] at(1.5,-1){\ };
	\end{tikzpicture}
	\quad
	\begin{tikzpicture}[scale=0.7]
	\draw[black, thick](2,3.46) arc (60:120:2);
	\draw[black, thick](0,0) arc (240:300:2);
	\draw[black, thick,dashed] (0,3.46) arc (120:240:2);
	\draw[black, thick,dashed] (2,0) arc (300:330:2);
	\draw[black, thick,dashed] (2.73,2.73) arc (390:420:2);
	\filldraw[black] (2.73,2.73) circle (.1);
	\filldraw[black] (2.73,0.73) circle (.1);
	\draw [black, very thick,->](1,3.73) --(1,-0.27);
	\node[right] at (1,1.73){$c'$};
	\node[above] at (1,3.73){$\operatorname{sgn}(c') = - \operatorname{sgn}(c)$};
	\node[below] at (1,-0.27){\,};
	\node[below] at (1,-1){$G(K')$};
\end{tikzpicture}
	\caption{The Gauss diagrams of $K$ and $K'$.}
		\label{fig18}
\end{figure}
Assume the chords corresponding to the crossing change of $K$ and $K'$ are $c$ and $c'$, respectively, $e$ is a chord of $K$ such that $e\ne c$, and $e'$ is the corresponding chord to $e$ in $K'$ after the crossing change. Then
$$
\operatorname{sgn}(c) = -\operatorname{sgn}(c'), \quad \operatorname{sgn}(e) = \operatorname{sgn}(e'), \quad d(c) = -d(c'), \quad d(e)=d(e'),
$$
since if $c$ intersects $e$, then $c'$ intersects $e$ in the opposite direction and has opposite sign. Moreover, since $\phi_c(k) = k \operatorname{mod} |d(c)|$ and $\phi_{c'} (k) = k \operatorname{mod} |d(c')|$ the property $d(c) = - d(c')$ implies that functions $\phi_c$ and $\phi_{c'}$ coincide.  For given $n$ denote $r^n (c) = \ell^n (c') = \{ \alpha_{1}, \ldots \alpha_{p_n} \}$ and $\ell^n (c) = r^n (c') = \{ \beta_1, \ldots, \beta_{q_n} \}$.
Then
$$
\operatorname{Ind}_{c}^{n} (z) =  \sum_{i=1}^{p_n} \operatorname{sgn} (\alpha_i) \, z^{\,  \phi_c (d(\alpha_i)) } - \sum_{j=1}^{q_n} \operatorname{sgn} (\beta_j) \, z^{\,  \phi_c (-d(\beta_j))}
$$
and
$$
\operatorname{Ind}_{c'}^{n} (z) =  \sum_{i=1}^{q_n} \operatorname{sgn} (\beta_i) \, z^{\,  \phi_{c'} (d(\beta_i)) } - \sum_{j=1}^{p_n} \operatorname{sgn} (\alpha_j) \, z^{\,  \phi_{c'} (-d(\alpha_j))},
$$
then
$$
\operatorname{Ind}_c^n(z) = - \operatorname{Ind}_{c'}^n(z^{-1}).
$$
By similar arguments we get
$$
\operatorname{Ind}_e^n(z) = \operatorname{Ind}_{e'}^n(z).
$$
Therefore,
\begin{alignat}{2}
H_K(t,y,z)-H_{K'}(t,y,z)
&=\sum_{n\in \mathbb{N}} \operatorname{sgn}(c) \left( t^{\operatorname{Ind}_c^n(z)} - 1\right) y^n - \sum_{n\in \mathbb{N}} \operatorname{sgn}(c') \left( t^{\operatorname{Ind}_{c'}^n(z)} - 1\right) y^n \nonumber\\
&=\sum_{n\in \mathbb{N}} \varepsilon \left(t^{\operatorname{Ind}_c^n(z)} + t^{- \operatorname{Ind}_c^n(z^{-1})}-2\right) y^n \nonumber,
\end{alignat}
where $\varepsilon = \operatorname{sgn} (c) = - \operatorname{sgn} (c')$. Thus, the Lemma is proved.
\end{proof}

By Lemma \ref{lem5.1-CCmove}, if $K$ and $K'$ differ by a crossing change move at a crossing, then the corresponding invariants will differ by a term whose coefficient is the sign of the chord corresponding to the crossing. This gives  a lower bound on the Gordian distance and proves Theorem~\ref{thm-dist}.

\begin{proof}
By the definition of homotopic knotoid diagrams, we may assume that $K'$ can be obtained by operating $k$ crossing change moves at $k$ crossings $c_1,\ldots, c_k$ on $K$.
Hence $d_G(K,K')\geq k$.
By Lemma \ref{lem5.1-CCmove}, we have 
$$
H_K(t,y,z) - H_{K'}(t,y,z) = \sum_{j=1}^{k} \left( \sum_{n \in \mathbb{N}} \operatorname{sgn}(c_j) \left( t^{\operatorname{Ind}^n_{c_j} (z)}+ t^{- \operatorname{Ind}^n_{c_j}(z^{-1})} - 2\right) y^n\right).
$$
Notice that there are only finite terms on the right hand side of this equality, hence we can change order of the two summing operations:
$$
H_K(t,y,z) - H_{K'}(t,y,z) = \sum_{n\in \mathbb{N}} \left(\sum_{j=1}^{k} \operatorname{sgn}(c_j) \left( t^{\operatorname{Ind}^n_{c_j}(z)} + t^{- \operatorname{Ind}^n_{c_j}(z^{-1})} - 2\right) \right)y^n.
$$
Some of the factors $\left(t^{z_n^m} + t^{-(z^{-1})_n^m}-2\right)$ might be zero, so we may rewrite the equality above as following:
$$
H_K(t,y,z) - H_{K'}(t,y,z)=\sum_{n\in \mathbb{N}} \left( \sum_{m\in\mathbb{N}} a_{n_m} \left( t^{z_n^m} + t^{-(z^{-1})_n^m}-2\right)\right)y^n,
$$
where $z_n^m$ indicates that the exponent of $t$ is a polynomial in $z$, the specific form of which depends on the given $n$ and $m$. We proved the first part of the theorem. For each $n\in\mathbb{N}$, there are no more than $k$ terms in which $\left | a_{n_m} \right | =1$, hence $d_G(K,K')\geq k \geq \sum\limits_{m \in\mathbb{N}}\left | a_{n_m} \right |$, which proves the second part of Theorem~\ref{thm-dist}. 
\end{proof}

\section{Examples} \label{sec4}

By point (3) of Theorem~\ref{thm-prop}, if $D$ is a zero height knotoid diagram then we have  $H_D(t,y,z) = 0$.  The following Example~\ref{ex:1} shows that there is a non-zero height knotoid diagram with $H_D(t,y,z) = 0$.

\begin{example} \label{ex:1} {\rm
 Let $D = 2_2$ be a diagram of knotoid $2_2$ from the Goundaroulis-Dorier-Stasiak list of planar knotoids~\cite{GDS19}. This diagram together with its  Gauss diagram, is shown  in Figure~\ref{fig19}. It is clear from the figure, that  $2_2$ is a non-zero height knotoid diagram.		
\begin{figure}[htbp]
\begin{center}
\tikzset{every picture/.style={line width=1pt}}
\scalebox{1.0}{
\begin{tikzpicture}[x=0.75pt,y=0.75pt,yscale=-0.75,xscale=0.75]
	\draw[black, very thick]   (184,794) .. controls (174,792) and (162,793) .. (161,773) .. controls (160,753) and (161,740) .. (171,733) .. controls (180.5,726) and (194,719) .. (230,720) .. controls (266,721) and (281,730) .. (284,759) .. controls (287,788) and (298,876) .. (270,876) .. controls (242,876) and (226,883) .. (224,856) ;
	\draw[black, very thick]   (191,825) .. controls (188,793) and (188,813) .. (188,785) .. controls (188,757) and (201,763) .. (216,764) .. controls (231,765) and (226,798) .. (225,839) ;
	\draw[black, very thick]  [shift={(191.16,825.63)}, rotate = 85.1]    (6,-3) .. controls (4,-1.4) and (2,-0.4) .. (0,0) .. controls (2,0.4) and (4,1.4) .. (6,3)   ;
	\draw[very thick]     (193,796) -- (204,797) -- (215,797) ;
	\filldraw[very thick]  (215,797) circle (2)   ;
	\draw[black, very thick]     (191,827) .. controls (194,861) and (215,843) .. (248,848) ;
	\filldraw[black, very thick]  [shift={(248,848)}, rotate = 8.62]  (0, 0) circle (2)  ;
	\filldraw[very thick]   (248,848) circle (2)  ;
	\draw (227,823) node [anchor=north west][inner sep=0.75pt]   [align=left] {$c_1$};
	\draw (168,798) node [anchor=north west][inner sep=0.75pt]   [align=left] {$c_2$};	
\end{tikzpicture}
\quad \quad \quad
\begin{tikzpicture}[scale=0.7]
	\draw[black, thick](0,0) arc (30:330:2);
	\filldraw[black] (0,0) circle (.1);
	\filldraw[black] (0,-2) circle (.1);
	\draw[black, very thick, ->] (-0.73,0.73) -- (-2.73,-2.73);
	\draw[black, very thick, <-] (-0.73,-2.73) -- (-2.73,0.73);	
	\node [above] at (-0.7,0.73){$+$};
	\node [above] at (-2.73,0.73){$-$};
	\node [left] at(-2.1,-0.42){$c_2$};
	\node [right] at(-1.3,-0.42){$c_1$};
\end{tikzpicture}
}
\caption{A knotoid diagram $D = 2_2$ and its Gauss diagram $G(D)$.}
\label{fig19}
\end{center}
\end{figure}
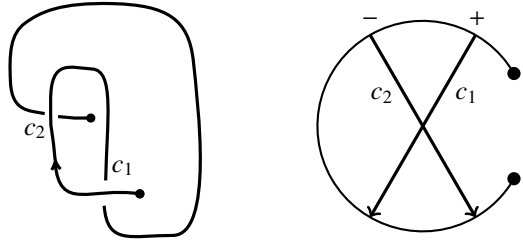
It is easy to see that $\operatorname{sgn} (c_1)=1$,  $\operatorname{sgn} (c_2)=-1$ and $d(c_1)=1$,   $d(c_2)=1$, hence $n = \gcd(d(c_1), d(c_2))=1$. Moreover, $\operatorname{Ind}_{c_1}^1(z)=1$ and  $\operatorname{Ind}_{c_2}^1(z)=1$. Therefore, $H_D(t,y,z) = \sum \operatorname{sgn} (c) (t^{\operatorname{Ind}^n_c (z)}-1)y^n = (t -1)y - (t-1)y = 0$. 
	}
\end{example}

The following Example~\ref{ex:2} demonstrates a calculation of an invariant~$H_D(t,y,z)$ according to (\ref{eqn:H}).

\begin{example}\label{ex:2} {\rm
	Let $D = D(5.1.28)$ be a diagram of  knotoid 5.1.28 from the Bartholomew's table~\cite{Bar15}. Diagram $D$ and its Gauss diagram $G(D)$ are presented in Figure~\ref{fig20}.
	\begin{figure}[htbp]
\begin{center}
	\tikzset{every picture/.style={line width=1pt}}	
	\begin{tikzpicture}[x=0.6pt,y=0.6pt,yscale=-1,xscale=1]
	\draw[black, very thick]    (54.41,2866.07) .. controls (82.28,2883.13) and (95.15,2877.22) .. (121.59,2848.35) ;
	\filldraw[black] (121.59,2848.35) circle (1.5pt);
	\draw[black, very thick] (129.05,2815.99) .. controls (148.35,2838.96) and (141.12,2863.03) .. (129.68,2880.09) ;
	\draw[black, very thick]    (125.8,2891.24) .. controls (110.08,2931.27) and (159.39,2964.73) .. (182.26,2956.86) .. controls (205.13,2948.98) and (210.13,2920.11) .. (207.27,2909.61) .. controls (204.41,2899.12) and (198.7,2885.34) .. (184.4,2885.99) .. controls (170.11,2886.65) and (156.53,2929.95) .. (146.53,2943.08) ;
	\draw[black, very thick]   (138.59,2949.93) .. controls (60.7,2978.14) and (31.4,2873.16) .. (55.69,2857.41) .. controls (79.99,2841.66) and (92.14,2847.57) .. (103.58,2860.69) ;
	\draw[black, very thick]    (109.37,2866.84) .. controls (117.94,2883.25) and (145.1,2900.96) .. (164.39,2903) ;
	\draw[black, very thick]    (172.25,2903.59) .. controls (182.97,2907.53) and (185.12,2906.87) .. (194.41,2904.24) ;
	\filldraw[black] (194.41,2904.24) circle (1.5pt);
	\draw[black, very thick]  (45.72,2856.49) .. controls (20.75,2818.51) and (84.98,2793.06) .. (127.76,2815.3) ;
	\draw[black, very thick, ->] (128,2815) -- (131,2820);
\node at (105,2880)  {$c_1$};
\node at (30,2860) {$c_2$};
\node at (140,2885)  {$c_3$};
\node at (140,2960) {$c_4$};
\node at (175,2915) {$c_5$};
\node [below]at (135,2975) {$D$};
	\end{tikzpicture}
\quad\quad
	\begin{tikzpicture}[scale=0.7]
		\draw[black, thick](0,0) arc (30:330:2);
		\filldraw[black] (0,0) circle (.1);
		\filldraw[black] (0,-2) circle (.1);
		\draw[black, very thick, ->] (-0.15,0.2) -- (-2.73,-2.73);
		\draw[black, very thick, ->] (-3.46,-2) -- (-0.73,0.73);
		\draw[black, very thick, ->] (-1.73,-3) -- (-1.73,1);
		\draw[black, very thick, ->] (-2.73,0.73) -- (-3.73,-1);
		\draw[black, very thick, ->] (-3.46,0) -- (-0.73,-2.73);
		\node [above]at(0,0){$+$};
		\node [above]at(-2.73,0.73){$-$};
		\node [above]at(-3.6,-0.1){$+$};
	    \node [below]at(-3.6,-1.8){$+$};
	    \node [below]at(-1.73,-3){$-$};
		\node [below]at(-1.5,-4){$G(D)$};
		\node [right]at(-0.8,-0.7){$c_1$};
		\node [right]at(-1.5,-0.2){$c_2$};
		\node [left]at(-1.7,0.3){$c_3$};
		\node [right]at(-3.7,-0.8){$c_4$};
		\node [below]at(-0.7,-1.8){$c_5$};		
	\end{tikzpicture}
	\caption{Knotoid diagram $D = D(5.1.28)$ and its Gauss diagram $G(D)$.}
	\label{fig20}
	\end{center}
\end{figure}
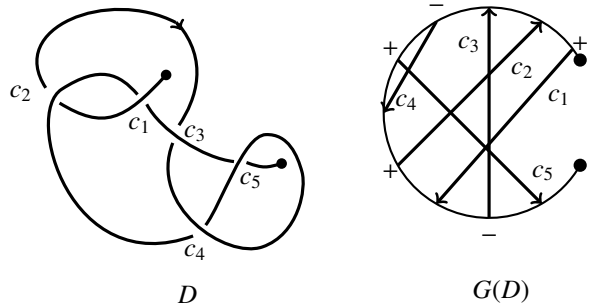
	One can see that $c_3$ and $c_5$ respectively pass through $c_1$ from left to right and from right to left along the direction of $c_1$, hence $d(c_1)=-2$. Analogously,  $c_3$ and $c_5$ respectively pass through $c_2$ from right to left and from left to right along the direction of $c_2$, hence $d(c_2)=2$. Similarly, it can be concluded that $d(c_3)=1$, $d(c_4)=-1$, and $d(c_5)=0$. Then we obtain
\begin{eqnarray*}
\operatorname{Ind}_{c_1}^1(z) & = & \operatorname{sgn}(c_3)\, z^{\phi _{c_1}(d(c_3))}=-z, \cr
\operatorname{Ind}_{c_1}^2(z) & = & -\operatorname{sgn}(c_5)\, z^{\phi _{c_1}(-d(c_5))}=-1, \cr
\operatorname{Ind}_{c_2}^1(z) & = & -\operatorname{sgn}(c_3)\, z^{\phi _{c_2}(-d(c_3))}=z^{-1}, \cr
\operatorname{Ind}_{c_2}^2(z) & = & \operatorname{sgn}(c_5)\, z^{\phi _{c_2}(d(c_5))}=1, \cr
\operatorname{Ind}_{c_3}^1(z) & = & \operatorname{sgn}(c_5) \, z^{\phi _{c_3}(d(c_5))} + \operatorname{sgn}(c_2) \, z^{\phi _{c_3}(d(c_2))}-\operatorname{sgn}(c_1) \, z^{\phi _{c_3}(-d(c_1))}=1, \cr
\operatorname{Ind}_{c_4}^1(z) & = & -\operatorname{sgn}(c_5) \, z^{\phi _{c_4}(-d(c_5))}=-1,\cr
\operatorname{Ind}_{c_5}^1(z) & = & \operatorname{sgn}(c_4) \, z^{\phi _{c_5}(d(c_4))} - \operatorname{sgn}(c_3) \, z^{\phi _{c_5}(-d(c_3))}=0, \cr
\operatorname{Ind}_{c_5}^2(z) & = & \operatorname{sgn}(c_1) \, z^{\phi _{c_5}(d(c_1))} - \operatorname{sgn}(c_2) \, z^{\phi _{c_5}(-d(c_2))}=0.
\end{eqnarray*}
Therefore
\begin{eqnarray*}
H_D(t,y,z) &=  & \sum_{\substack{c\in C(G\left ( D \right ))\\n\in \mathbb{N}}}\operatorname{sgn} ( c ) \left( t^{\operatorname{Ind}_{c}^{n}  ( z )} - 1   \right ) \, y^n \cr
		& = & \left( t^{-z}-1\right)\, y \, + \, \left( t^{-1} - 1 \right) \, y^2 \, + \, \left( t^{z^{-1}}-1 \right) \, y \, + \, (x-1) \, y^2 \, \cr
		& &
		- \, (t-1) \, y \, - \, (t^{-1}-1) \, y \cr
		&= & \left[ \left( t^{-z} + t^{z^{-1} } \, - \, 2\right) \, - \, \left( t + t^{-1} - 2\right) \right] \, y \, + \, \left( t + t^{-1} - 2\right) y^2.
\end{eqnarray*}		
	}
\end{example}

The following Example~\ref{ex:3} demonstrates that $H_D(t,y,z)$ distinguishes a knotoid and its inverse image.

\begin{example}\label{ex:3} {\rm
	Let us consider the knotoid 5.1.28 of Bartholomew's table~\cite{Bar15} and its inverse. A diagram $D = D(5.1.28)$ of the knotoid 5.1.28 is presented in Figure~\ref{fig20}. It was shown in Example~\ref{ex:2} that  	
$$
H_D(t,y,z) = \left[ \left( t^{z^{-1}} + t^{-z}-2\right) - \left( t+t^{-1}-2 \right) \right] \, y \, + \, \left( t+t^{-1}-2 \right) \, y^2.
$$

The inverse diagram $-D = -D(5.1.28)$ and its Gauss diagram $G(-D)$ are presented in Figure~\ref{fig21}.
	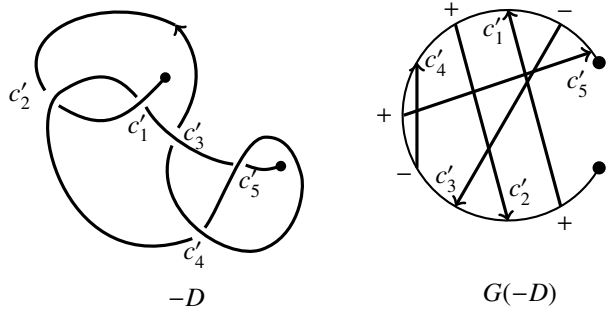
\begin{figure}[htbp]
	\begin{center}
		\tikzset{every picture/.style={line width=1pt}}	
		\begin{tikzpicture}[x=0.6pt,y=0.6pt,yscale=-1,xscale=1]
			\draw[black, very thick]    (54.41,2866.07) .. controls (82.28,2883.13) and (95.15,2877.22) .. (121.59,2848.35) ;
			\filldraw[black] (121.59,2848.35) circle (1.5pt);
			\draw[black, very thick] (129.05,2815.99) .. controls (148.35,2838.96) and (141.12,2863.03) .. (129.68,2880.09) ;
			\draw[black, very thick]    (125.8,2891.24) .. controls (110.08,2931.27) and (159.39,2964.73) .. (182.26,2956.86) .. controls (205.13,2948.98) and (210.13,2920.11) .. (207.27,2909.61) .. controls (204.41,2899.12) and (198.7,2885.34) .. (184.4,2885.99) .. controls (170.11,2886.65) and (156.53,2929.95) .. (146.53,2943.08) ;
			\draw[black, very thick]   (138.59,2949.93) .. controls (60.7,2978.14) and (31.4,2873.16) .. (55.69,2857.41) .. controls (79.99,2841.66) and (92.14,2847.57) .. (103.58,2860.69) ;
			\draw[black, very thick]    (109.37,2866.84) .. controls (117.94,2883.25) and (145.1,2900.96) .. (164.39,2903) ;
			\draw[black, very thick]    (172.25,2903.59) .. controls (182.97,2907.53) and (185.12,2906.87) .. (194.41,2904.24) ;
			\filldraw[black] (194.41,2904.24) circle (1.5pt);
			\draw[black, very thick]  (45.72,2856.49) .. controls (20.75,2818.51) and (84.98,2793.06) .. (127.76,2815.3) ;
			\draw[black, very thick, <-] (128,2815) -- (131,2819);
			\node at (105,2880)  {$c'_1$};
			\node at (30,2860) {$c'_2$};
			\node at (140,2885)  {$c'_3$};
			\node at (140,2960) {$c'_4$};
			\node at (175,2915) {$c'_5$};
			\node [below]at (135,2975) {$-D$};
		\end{tikzpicture}
		\quad\quad
		\begin{tikzpicture}[scale=0.7]
			\draw[black, thick](0,0) arc (30:330:2);
			\filldraw[black] (0,0) circle (.1);
			\filldraw[black] (0,-2) circle (.1);
			\draw[black, very thick, ->] (-3.71,-1) -- (-0.15,0.2);
			\draw[black, very thick, ->] (-0.73,0.73) -- (-2.73,-2.73);
			\draw[black, very thick, ->] (-0.73,-2.73) -- (-1.73,1);
			\draw[black, very thick, ->] (-2.73,0.73) -- (-1.73,-3);
			\draw[black, very thick, <-] (-3.46,0) -- (-3.46,-2);
			\node [above]at(-0.65,0.65){$-$};
			\node [above]at(-2.8,0.7){$+$};
			\node [left]at(-3.73,-1){$+$};
			\node [below]at(-3.7,-1.8){$-$};
			\node [below]at(-0.65,-2.73){$+$};
			\node [below]at(-1.5,-4){$G(-D)$};
			\node [right]at(-2.4,0.6){$c'_1$};
			\node [below]at(-1.5,-2){$c'_2$};
			\node [right]at(-3.3,-2.2){$c'_3$};
			\node [right]at(-3.5,0){$c'_4$};
			\node [left]at(0,-0.4){$c'_5$};		
		\end{tikzpicture}
		\caption{Knotoid diagram $-D$ and its Gauss diagram $G(-D)$.} \label{fig21}
	\end{center}
\end{figure}

By the direct  calculations, we get for this diagram
$$
d(c'_1)=2, \quad d(c'_2)=-2, \quad d(c'_3)=-1, \quad d(c'_4)=1, \quad d(c'_5)=0.
$$
Moreover,
$$
\operatorname{Ind}_{c'_1}^1(z)=z, \quad \operatorname{Ind}_{c'_1}^2(z)=1, \quad \operatorname{Ind}_{c'_2}^1(z)=-z^{-1}, \quad \operatorname{Ind}_{c'_2}^2(z)=-1,
$$
and
$$
\operatorname{Ind}_{c'_3}^1(z)=-1, \quad \operatorname{Ind}_{c'_4}^1(z)=1, \quad \operatorname{Ind}_{c'_5}^1(z)=0, \quad \operatorname{Ind}_{c'_5}^2(z)=0.
$$
Therefore,
\begin{eqnarray*}
H_{-D}(t,y,z) &= & ( t^{z}-1 ) \, y \, + \, ( t-1 ) \, y^2 \, + \, \left( t^{-z^{-1}}-1 \right) \, y \, +\, \left( t^{-1}-1 \right) y^2\, \cr
 & & - \, \left( t^{-1}-1\right)\, y \, - \, ( t-1) \, y \cr
	& = & \left[ \left( t^{z}+t^{-z^{-1}}-2 \right) \, - \, \left( t+t^{-1}-2 \right) \right] \, y \, + \, \left( t+t^{-1}-2 \right) \, y^2.
\end{eqnarray*}
Since $H_{D}(t,y,z) \ne H_{-D}(t,y,z)$, we conclude that $D = D(5.1.28)$ is irreversible.
}
\end{example}

In Example~\ref{exam4.3}, we calculate that
$H_{K(c)}(t,y,z)\ne 0$ for a singular planar knotoid $K(c)$ with one singular crossing $c$. Thus $H_K(t,y,z)$ is a Vassiliev invariant of order one for planar knotoids.

\begin{example}\label{exam4.3} {\rm
Let $K(c)$ be a singular planar knotoid with only one singular crossing as shown in Figure~\ref{fig22}. For $K(c^+)$ the direct calculations give
$$
d(c_1)=-2, \quad d(c_2)=-2, \quad d(c_3)=1, \quad d(c_4)=1,
$$
and
$$
\begin{gathered}
\operatorname{Ind}_{c_1}^1(z)=-z, \quad \operatorname{Ind}_{c_1}^2(z)=-1, \quad \operatorname{Ind}_{c_2}^1(z)=-z, \cr  \operatorname{Ind}_{c_2}^2(z)=-1, \quad \operatorname{Ind}_{c_3}^1(z)=1, \quad \operatorname{Ind}_{c_4}^1(z)=1.
\end{gathered}
$$
Hence $H_{K(c^+)}(t,y,z)=0$.

For $K(c^-)$, the direct calculations give
$$
d(c_1)=-2, \quad d(c_2)=2, \quad d(c_3)=1, \quad d(c_4)=1
$$
and
$$
\begin{gathered}
\operatorname{Ind}_{c_1}^1(z)=-z, \quad \operatorname{Ind}_{c_1}^2(z)=-1, \quad \operatorname{Ind}_{c_2}^1(z)=z, \cr  \operatorname{Ind}_{c_2}^2(z)=1, \quad \operatorname{Ind}_{c_3}^1(z)=1, \quad \operatorname{Ind}_{c_4}^1(z)=1.
\end{gathered}
$$
Hence $H_{K(c^-)}(t,y,z) = -(t^{-z}-1)y-(t^{-1}-1)y^2-(t^z-1)y-(t-1)y^2$.

According to the Vassiliev skein relation (\ref{eqn:Vas}), we calculate that
\begin{eqnarray*}
H_{K(c)}(t,y,z) & = & H_{K(c^+)}(t,y,z) \, - \, H_{K(c^-)}(t,y,z) \, \cr
& = &  (t^{-z}-1)y+(t^{-1}-1)y^2+(t^z-1)y+(t-1)y^2.
\end{eqnarray*}
\begin{figure}[htbp]
	\begin{center}
		\begin{tikzpicture}[x=0.75pt,y=0.75pt,yscale=-0.5,xscale=0.5]
			\draw[black, very thick] (256.4,316.05) .. controls (257.3,301.72) and (258.79,230.07) .. (302.71,265.35) ;
			\filldraw[black] (256.4,316.05) circle (3pt);
			\draw[black, very thick] (302.71,265.35) -- (352,317)  ;
			\draw[black, very thick] (360.07,324.87) .. controls (400.7,359.04) and (456.27,311.64) .. (416.84,273.06);
			\draw[black, very thick] (255.51,271.96) .. controls (213.38,265.35) and (226.83,361.24) .. (287.77,378.88) .. controls (348.72,396.52) and (365.75,245.5) .. (409,267) ;
			\draw[black, very thick] (268.26,273.96) .. controls (295.33,280.65) and (298.71,269.35) .. (302.71,265.35) ;
			\draw[black, very thick] (363.33,226.65) .. controls (404.33,218.65) and (429.33,258.65) .. (403.39,292.9) ;
			\filldraw[black] (403.39,292.9) circle (3pt);
			\draw[black, very thick] (302.71,265.35) .. controls (323.4,241.71) and (331.1,238.74) .. (361.91,227.18) ;
			\draw [shift={(363.33,226.65)}, rotate = 159.44] [black, very thick] (6.56,-2.94) .. controls (4.17,-1.38) and (1.99,-0.4) .. (0,0) .. controls (1.99,0.4) and (4.17,1.38) .. (6.56,2.94)   ;
			\filldraw[black] (302,265.35) circle (4pt);
			\node at (313,405) {$K(c)$};
			\node at (292,239) {$c$};
		\end{tikzpicture}
		\quad
		\begin{tikzpicture}[x=0.75pt,y=0.75pt,yscale=-0.5,xscale=0.5]
			\draw[black, very thick] (256.4,316.05) .. controls (257.3,301.72) and (258.79,230.07) .. (302.71,265.35) ;
			\filldraw[black] (256.4,316.05) circle (3pt);
			\draw[black, very thick] (302.71,265.35) -- (352,317)  ;
			\draw[black, very thick] (360.07,324.87) .. controls (400.7,359.04) and (456.27,311.64) .. (416.84,273.06);
			\draw[black, very thick] (255.51,271.96) .. controls (213.38,265.35) and (226.83,361.24) .. (287.77,378.88) .. controls (348.72,396.52) and (365.75,245.5) .. (409,267) ;
			\draw[black, very thick] (268.26,273.96) .. controls (295.33,280.65) and (298.71,269.35) .. (299,270) ;
			\draw[black, very thick] (363.33,226.65) .. controls (404.33,218.65) and (429.33,258.65) .. (403.39,292.9) ;
			\filldraw[black] (403.39,292.9) circle (3pt);
			\draw[black, very thick] (306,260) .. controls (323.4,241.71) and (331.1,238.74) .. (361.91,227.18) ;
			\draw [shift={(363.33,226.65)}, rotate = 159.44] [black, very thick] (6.56,-2.94) .. controls (4.17,-1.38) and (1.99,-0.4) .. (0,0) .. controls (1.99,0.4) and (4.17,1.38) .. (6.56,2.94)   ;
			\node at (313,405) {$K(c^+)$};
			\node at (292,239) {$c$};
		\end{tikzpicture}
		\quad
		\begin{tikzpicture}[x=0.75pt,y=0.75pt,yscale=-0.5,xscale=0.5]
			\draw[black, very thick] (256.4,316.05) .. controls (257.3,301.72) and (258.79,230.07) .. (300,260) ;
			\filldraw[black] (256.4,316.05) circle (3pt);
			\draw[black, very thick] (306,270) -- (352,317) ;
			\draw[black, very thick] (360.07,324.87) .. controls (400.7,359.04) and (456.27,311.64) .. (416.84,273.06);
			\draw[black, very thick] (255.51,271.96) .. controls (213.38,265.35) and (226.83,361.24) .. (287.77,378.88) .. controls (348.72,396.52) and (365.75,245.5) .. (409,267) ;
			\draw[black, very thick] (268.26,273.96) .. controls (295.33,280.65) and (298.71,269.35) .. (302.71,265.35) ;
			\draw[black, very thick] (363.33,226.65) .. controls (404.33,218.65) and (429.33,258.65) .. (403.39,292.9) ;
			\filldraw[black] (403.39,292.9) circle (3pt);
			\draw[black, very thick] (302.71,265.35) .. controls (323.4,241.71) and (331.1,238.74) .. (361.91,227.18) ;
			\draw [shift={(363.33,226.65)}, rotate = 159.44] [black, very thick] (6.56,-2.94) .. controls (4.17,-1.38) and (1.99,-0.4) .. (0,0) .. controls (1.99,0.4) and (4.17,1.38) .. (6.56,2.94)   ;
			\node at (313,405) {$K(c^-)$};
			\node at (292,239) {$c$};
		\end{tikzpicture}
		\caption{A singular planar knotoid $K(c)$ with one singular crossing and its Vassiliev resolutions $K(c^+)$ and $K(c^-)$.}
		\label{fig22}
	\end{center}
\end{figure}
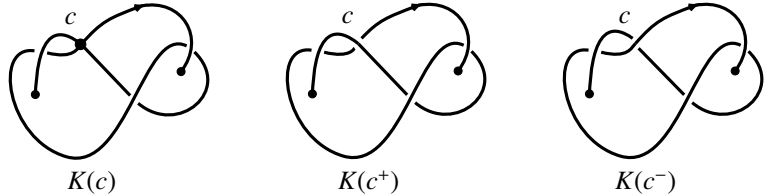

}
\end{example}

\begin{example} {\rm
Let $K$ be the planar knotoid 5.1.28 from the Bartholomew's table~\cite{Bar15}. Its Gauss diagram $D$ is shown in Figure~\ref{fig20}. It was calculated in Example~\ref{ex:2} that  	
$$
H_D(t,y,z) = [(t^{z} + t^{-z^{-1}}-2)-(t+t^{-1}-2)]y+(t+t^{-1}-2)y^2.
$$
Compare this expression with $H_{\bf 0}(t,y,z)=0$, where $\bf 0$ is a trivial planar knotoid. By Theorem~\ref{thm-dist} we get that $d_G(K,\bf 0) \geq 2$. At the same time, Figure~\ref{fig23} shows that $d_G(K, \bf 0)\le2$. Hence $d_G(K,\bf 0)=2$.
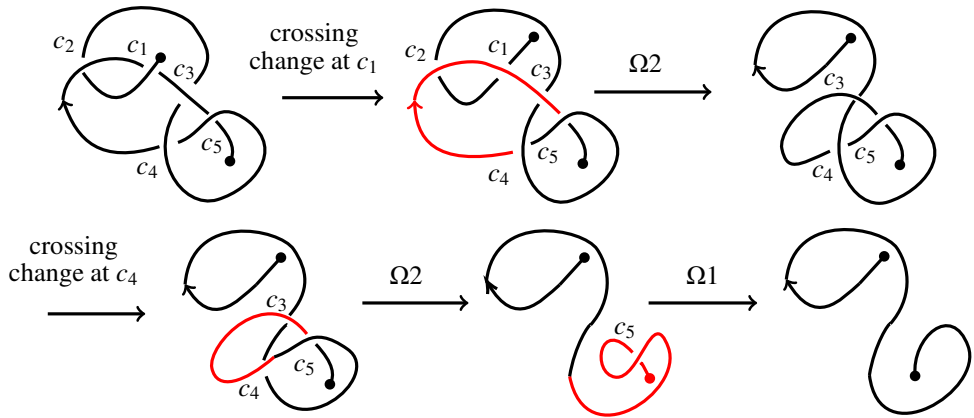
\begin{figure}[htbp]
\begin{center}
\tikzset{every picture/.style={line width=1pt}}
\begin{tikzpicture}[x=0.7pt,y=0.7pt,yscale=-1,xscale=1]
	\draw[black, very thick]    (55.08,82.34) .. controls (59.54,44.49) and (107.25,54.09) .. (114.83,65.28) .. controls (122.41,76.48) and (126.93,83.12) .. (115.33,96.45) ;
	\draw[black, very thick]    (55.08,90.51) .. controls (72.03,114.33) and (81.84,103.13) .. (97,82.34) ;
	\filldraw[black] (97,82.34)  circle (1.5pt);
	\draw[black, very thick]  (44.19,105.22) .. controls (45.72,85.01) and (72.47,77.28) .. (87.63,87.14) ;
	\draw[black, very thick]    (107.33,107.45) .. controls (87.71,134.64) and (108.59,164.97) .. (127.76,159.11) .. controls (146.94,153.24) and (164.33,132.45) .. (144.71,116.99) .. controls (125.09,101.53) and (121.52,123.39) .. (106.36,128.19) ;
	\draw[black, very thick]    (96.11,90.34) .. controls (106.81,98.33) and (104.13,97.09) .. (122.86,112.73) ;
	\draw[black, very thick]    (126.87,116.63) .. controls (136.68,127.12) and (137.13,133.52) .. (134.45,138.32) ;
	\filldraw[black] (134.45,138.32) circle (1.5pt);
	\draw[black, very thick, <-] (44.19,105.22) .. controls (51.87,127.4) and (64.32,140.15) .. (96.33,130.45) ;
	\draw[black,very thick, ->]    (162.67,104) -- (217.33,104) ;
	\draw [red, very thick, <-] (234.2,105.36) .. controls (234.94,125.67) and (243.91,142.13) .. (287.33,133.45) ;
	\draw[black, very thick]  (245.75,83.68) .. controls (250.21,45.83) and (297.92,55.42) .. (305.5,66.62) .. controls (313.08,77.81) and (316.93,87.12) .. (305.33,100.45) ;
	\draw[black, very thick]    (247.33,92.45) .. controls (262.25,119.72) and (268.67,103.11) .. (279.33,93.78) ;
	\draw [red, very thick]   (234.2,105.36) .. controls (236.39,85.01) and (265.33,82.45) .. (273.33,85.78) .. controls (281.33,89.11) and (288,88.45) .. (312.67,112.45) ;
	\draw[black, very thick]  (301.33,107.45) .. controls (281.71,134.64) and (299.26,165.64) .. (318.43,159.77) .. controls (337.6,153.91) and (354.99,133.12) .. (335.37,117.66) .. controls (315.75,102.2) and (312.19,124.06) .. (297.03,128.85);
	\draw[black, very thick]  (317.54,117.97) .. controls (327.35,128.46) and (327.79,134.85) .. (325.12,139.65) ;
	\filldraw[black] (325.12,139.65) circle (1.5pt);
	\draw[black, very thick] (298.67,71.11) -- (285.33,85.78) ;
	\filldraw[black] (298.67,71.11) circle (1.5pt);
	\draw[black, very thick] (418.2,86.31) .. controls (421.54,46.49) and (469.25,56.09) .. (476.83,67.28) .. controls (484.41,78.48) and (486.2,89.67) .. (474.6,103) ;
	\draw[black, very thick, <-] (418.2,86.31) .. controls (432.7,110.82) and (446.27,95.55) .. (456.67,86.45) ;
	\draw [black, very thick]   (458.67,132.89) .. controls (441.33,140.89) and (435.33,144.23) .. (432.67,131.56) .. controls (430,118.89) and (462,84.89) .. (484,113.11);
	\draw[black, very thick]  (471.93,105.53) .. controls (452.31,132.72) and (470.59,166.3) .. (489.76,160.44) .. controls (508.94,154.58) and (526.33,133.79) .. (506.71,118.33) .. controls (487.09,102.87) and (483.52,124.72) .. (468.36,129.52) ;
	\draw[black, very thick] (488.87,118.63) .. controls (498.68,129.12) and (499.13,135.52) .. (496.45,140.32) ;
	\filldraw[black] (496.45,140.32) circle (1.5pt);
	\draw[black, very thick] (470,71.78) -- (456.67,86.45) ;
	\filldraw [black] (470,71.78) circle (1.5pt);
	\draw[black, very thick, ->] (331.33,100.89) -- (394.67,100.89) ;
	\draw[black,very thick, ->] (36,220.89) -- (88.67,220.89) ;
	\draw[black, very thick] (110.2,204.97) .. controls (113.54,165.16) and (161.25,174.75) .. (168.83,185.95) .. controls (176.41,197.14) and (178.2,208.34) .. (166.6,221.67) ;
	\draw[black, very thick,<-]    (110.2,204.97) .. controls (124.7,229.49) and (138.27,214.21) .. (148.67,205.11);
	\draw [red, very thick]   (158.36,244.19) .. controls (142.36,258.19) and (127.33,262.89) .. (124.67,250.23) .. controls (122,237.56) and (154,203.56) .. (176,231.78) ;
	\draw[black,very thick] (154,254.67) .. controls (168.67,296) and (218.67,254.67) .. (198,240) .. controls (177.33,225.34) and (173.52,239.39) .. (158.36,244.19) ;
	\draw[black, very thick]  (180.87,237.3) .. controls (190.68,247.79) and (191.13,254.19) .. (188.45,258.98) ;
	\filldraw[black] (188.45,258.98) circle (1.5pt);
	\draw[black, very thick] (162,190.45) -- (148.67,205.11) ;
	\filldraw[black](162,190.45) circle (1.5pt);
	\draw[black, very thick] (165.33,226.89) .. controls (164.67,226) and (156.33,235.78) .. (152.67,246) ;
	\draw[black, very thick, ->] (206,215) -- (263.33,215);
	\draw[black, very thick] (273.08,203.01) .. controls (277.54,165.16) and (325.25,174.75) .. (332.83,185.95) .. controls (340.41,197.14) and (340.93,213.57) .. (329.33,226.89) ;
	\draw [black, very thick] (274.2,204.97) .. controls (288.7,229.49) and (302.27,214.21) .. (312.67,205.11) ;
	\draw [shift={(273.08,203.01)}, rotate = 61.32] [black, very thick] (6.56,-2.94) .. controls (4.17,-1.38) and (1.99,-0.4) .. (0,0) .. controls (1.99,0.4) and (4.17,1.38) .. (6.56,2.94)   ;
	\draw [red, very thick]   (318,254.67) .. controls (326.67,298.89) and (380,263.11) .. (371.33,239.11) .. controls (362.67,215.11) and (356,257.11) .. (342.67,255.11) .. controls (329.33,253.11) and (334,221.78) .. (352,241.78) ;
	\draw [red, very thick  ,draw opacity=1 ]   (356.21,247.3) .. controls (358,249.11) and (355.33,244.45) .. (361.33,255.78) ;
	\filldraw[red] (361.33,255.78) circle (1.5pt);
	\draw[black, very thick]    (326,190.45) -- (312.67,205.11) ;
	\filldraw[black] (326,190.45) circle (1.5pt);
	\draw[black, very thick]    (329.33,226.89) .. controls (328.67,226) and (320.67,238) .. (318,254.67) ;
	\draw[black, very thick, ->]    (360.67,215) -- (419.33,215) ;
	\draw[black, very thick]  (435.08,202.34) .. controls (439.54,164.49) and (487.25,174.09) .. (494.83,185.28) .. controls (502.41,196.48) and (502.93,212.9) .. (491.33,226.23);
	\draw[black, very thick, <-] (435.08,202.34) .. controls (450.7,228.82) and (464.27,213.55) .. (474.67,204.45) ;
	\draw[black, very thick]   (480,254) .. controls (485.33,297.56) and (542,262.45) .. (533.33,238.45) .. controls (524.67,214.45) and (502.67,233.56) .. (504.67,254.45) ;
	\filldraw[black] (504.67,254.45) circle (1.5pt);
	\draw[black, very thick] (488,189.78) -- (474.67,204.45) ;
	\filldraw[black] (488,189.78) circle (1.5pt);
	\draw[black, very thick] (491.33,226.23) .. controls (490.67,225.34) and (482.67,237.34) .. (480,254) ;	
	\node at (85,77) {$c_1$};
	\node at  (45,75) {$c_2$};
	\node at (110,90) {$c_3$};
	\node at (90,143) {$c_4$};
	\node at (125,130) {$c_5$};
	\node at (180,70) {crossing};
	\node at (180,85) {change at $c_1$};
	\node at (280,77) {$c_1$};
	\node at (235,80) {$c_2$};
	\node at (303,90) {$c_3$};
	\node at (280,147) {$c_4$};
	\node at (306.67,135) {$c_5$};
	\node at (461.33,95) {$c_3$};
	\node at (455,145){$c_4$};
	\node at (478,135) {$c_5$};
	\node at (359,85) {$\Omega2$};
	\node at (50,185) {crossing};
	\node at (50,200) {change at $c_4$};
	\node at (160,215) {$c_3$};
	\node at (145,263) {$c_4$};
	\node at (175,250) {$c_5$};
	\node at (229.67,200) {$\Omega2$};
	\node at (348,230) {$c_5$};
	\node at (390,200) {$\Omega1$};
\end{tikzpicture}
		\caption{A sequence of crossing changes and Reidemeister moves from $D(5.1.28)$ to a trivial planar knotoid.} \label{fig23}
	\end{center}
\end{figure}
}
\end{example}

\end{document}